\documentclass{article}
\usepackage{amsmath}
\usepackage{amssymb}
\usepackage{amsthm}
\usepackage[dvips]{graphics,color}

\theoremstyle{plain}
\newtheorem{Thm}{Theorem}
\newtheorem{Lem}{Lemma}[section]
\newtheorem{Prop}[Lem]{Proposition}
\newtheorem{Cor}[Lem]{Corollary}

\theoremstyle{definition}
\newtheorem{Def}[Lem]{Definition}
\newtheorem{Rem}[Lem]{Remark}

\theoremstyle{remark}
\newtheorem*{Rem*}{Remark}

\newtheorem*{Acknowledgments}{Acknowledgments}

\newcommand{\bbC}{\mathbb{C}}

\newcommand{\bbN}{\mathbb{N}}

\newcommand{\bbQ}{\mathbb{Q}}
\newcommand{\bbR}{\mathbb{R}}
\newcommand{\bbZ}{\mathbb{Z}}

\newcommand{\tW}{\Tilde{W}}

\newcommand{\cI}{\mathcal{I}}
\newcommand{\cM}{\mathcal{M}}

\newcommand{\cU}{\mathcal{U}}

\DeclareMathOperator{\area}{area}
\DeclareMathOperator{\dens}{dens}

\DeclareMathOperator{\im}{im}

\DeclareMathOperator{\sgn}{sgn}
\DeclareMathOperator{\spec}{spec}

\DeclareMathOperator{\Spec}{Spec}

\renewcommand{\leq}{\leqslant}
\renewcommand{\geq}{\geqslant}

\newcommand{\eps}{\varepsilon}
\newcommand{\vhi}{\varphi}

\renewcommand{\vec}[1]{{\langle{#1}\rangle}}

\renewcommand{\include}{\input}

\begin{document}
\title{Slowly divergent geodesics in moduli space}
\author{Yitwah Cheung}
\maketitle

\begin{abstract}\noindent
Slowly divergent Teichm\"uller geodesics in the moduli space of 
Riemann surfaces of genus $g\geq2$ are constructed via cyclic 
branched covers of the torus.  Nonergodic examples (i.e. geodesics 
whose defining quadratic differential has nonergodic vertical 
foliation) diverging to infinity at sublinear rates are 
constructed using a Diophantine condition.  Examples with an 
arbitrarily slow prescribed rate of divergence are also exhibited.  
\end{abstract}

\section{Introduction}
Let $\cM_g$ denote the moduli space of closed Riemann surfaces 
of genus $g\geq2$, endowed with the Teichm\"uller metric $\tau$.  
A geodesic in $\cM_g$ is determined by a pair $(X_0,q)$ where 
$X_0$ is a Riemann surface and $q$ is a holomorphic quadratic 
differential on $X_0$.  The differential $q$ defines a flat 
metric with isolated singularities on $X_0$ together with a 
pair of transverse measured foliations defined by $q>0$ (the 
horizontal) and $q<0$ (the vertical).  By a theorem of Masur 
\cite{Ma92} the vertical foliation of $q$ is uniquely ergodic 
if $X_t$ accumulates in $\cM_g$ as $t\to\infty$.  Therefore, 
a \textit{nonergodic} geodesic, by which we mean a geodesic 
determined by a pair $(X_0,q)$ such that the vertical foliation 
of $q$ is not uniquely ergodic, must eventually leave every 
compact set.  A geodesic with this latter property is said to 
be \textit{divergent}.  The original motivation of this study 
is to answer a question of C.~McMullen regarding the existence 
of \textit{slowly divergent} nonergodic geodesics: 
\begin{equation}\label{D:Roth}
\lim_{t\to\infty} \frac{\tau(X_0,X_t)}{t} = 0.  
\end{equation}

The examples are realized using branched covers of the torus 
satisfying a Diophantine condition.  Let $(X,q)$ be the $g$-cyclic 
branched cover of $T=(\bbC/\bbZ[i],dz^2)$ obtained by cutting 
along an embedded linear arc $\gamma$.  (See \S\ref{S:Cyclic} 
for a precise definition.)  Each $\theta\in S^1$ determines 
a Teichm\"uller geodesic $X_t^\theta$ in $\cM_g$ starting at 
$X_0^\theta=X$.  A direction $\theta$ is also said to be 
\textit{nonergodic, divergent, slowly divergent, etc.} if 
the corresponding geodesic $X_t^\theta$ has the same property.  
For a slowly divergent direction it makes sense to consider 
the \textit{sublinear rate}: 
\begin{equation}\label{D:Sublinear}
  r_+(\theta) := \limsup_{t\to\infty} 
                 \frac{\log\tau(X_t^\theta,X_0^\theta)}{\log t}
\end{equation}

A pair $(x_0,y_0)\in\bbR^2$ is said to satisfy a Diophantine 
condition if there are constants $c_0>0$ and $d_0>0$ such that 
for all pairs of integers $(m,n)\in\bbZ^2\setminus\{0\}$ 
\begin{equation}\label{DIO}
  \inf_{l\in\bbZ}|mx_0+ny_0+l|>\frac{c_0}{\max(|m|,|n|)^{d_0}}. 
\end{equation}

Let $x_0+iy_0\in\bbC$ be the affine holonomy $\int_\gamma dz$ of $\gamma$.  
\begin{Thm}\label{Thm1}
If $(x_0,y_0)$ satisfies a Diophantine condition with exponent $d_0$ 
then for every $e_0>\max(d_0,2)$ there is a Hausdorff dimension $1/2$ 
set of slowly divergent nonergodic directions $\theta$ with sublinear 
rate $r_+(\theta)\leq1-1/e_0$.  
\end{Thm}

It should be emphasized that Theorem~1 does not give any examples 
of nonergodic directions with $r_+(\theta)\leq1/2$.  In fact, after 
this paper had been accepted, it was shown that if $r_+(\theta)\leq1/2$ 
then $\theta$ is uniquely ergodic.  See \cite{CE}.  

As a complement to Theorem~1, we also prove 
\begin{Thm}\label{Thm2}
If $(x_0,y_0)\not\in\bbQ^2$ then there are directions which are 
divergent with an arbitrarily slow prescribed rate, i.e. given 
any function $R(t)$ with $R(t)\to\infty$ as $t\to\infty$ there 
exists a divergent direction $\theta$ such that 
$\tau(X_t^\theta,X_0^\theta)\leq R(t)$ for all sufficiently large $t$.   
\end{Thm}

In $\cM_1$, the asymptotic behavior of a geodesic is determined by 
the arithmetic properties of its endpoint in $\bbR\cup\{\infty\}$.  
For example, (\ref{D:Roth}) holds iff the endpoint is a Roth number: 
for any $\eps>0$, there exists $c_0>0$ such that 
for any $p,q\in\bbZ$, $|\alpha-p/q|>c_0/|q|^{2+\eps}$.  
The question asked by C.~McMullen was inspired by a recent result 
of Marmi-Moussa-Yoccoz concerning interval exchange maps, which 
give another a source of Teichm\"uller geodesics in $\cM_g$.  
In \cite{MMY}, they construct examples of uniquely ergodic interval 
exchange maps based on a certain ``Roth type'' condition, which is 
apparently stronger than (\ref{D:Roth}).  Theorem~1 shows that the 
condition (\ref{D:Roth}) alone is not sufficient to ensure unique 
ergodicity.  

\begin{Acknowledgments}
The author would like to thank the referees for numerous suggestions 
that significantly improved the exposition of the paper.  In addition 
thanks must go to Howard~Masur, Curt~McMullen and Barak~Weiss for many 
enlightening discussions on the topic of slow divergence.  Last, but 
not least, the author is indebted to his wife Ying~Xu for her 
constant and unwavering support.  
\end{Acknowledgments}

\section{Cyclic branched covers along a slit}\label{S:Cyclic}
The $g$-cyclic branched cover of $T$ along $\gamma$ is defined 
as follows.  Endow the complement of $\gamma$ with the metric 
defined by shortest path and let $T'$ be its metric completion.  
$T'$ is a compact surface with a single boundary component and 
is known as a \textit{slit torus}, i.e. $T$ slit along $\gamma$.  
Let $\gamma_\pm$ denote the two lifts of $\gamma$ under the 
natural projection $T'\to T$ which maps $\partial T'$ onto $\gamma$.  
For convenience, we assume $\gamma$ is defined on the unit interval.  
Let $X$ be the quotient space of $T'\times\bbZ/g\bbZ$ obtained 
by identifying $(\gamma_-(t),n)$ with $(\gamma_+(t),n+1)$ for 
all $t\in[0,1]$ and $n\in\bbZ/g\bbZ$.  
The map $\pi:X\to T$ induced by projection onto the first factor 
is a branched cover of degree $g$, holomorphic with respect to 
a unique complex structure on $X$.  The pair $(X,\pi^*dz^2)$ is 
called the $g$-cyclic cover of $T=(\bbC/\bbZ[i],dz^2)$ along $\gamma$.  

Note that $X$ is a closed Riemann surface of genus $g$.  
The map $\pi$ is branched at two points corresponding to 
the zeros of the quadratic differential $\pi^*dz^2$.  
Each branch point lies over an endpoint of $\gamma$.  

\subsection{Teichm\"uller geodesics and saddle connections}
The Teichm\"uller geodesic $X_t^\theta$ will be described explicitly.  
Let $g_t^\theta:\bbR^2\to\bbR^2$ be the linear map which contracts 
distances by a factor of $e^{t/2}$ in the $\theta$ direction while 
expanding by $e^{t/2}$ in the direction perpendicular to $\theta$.  
There is an atlas of charts $\{U_\alpha,\vhi_\alpha\}$ covering $X$ 
away from the branch points such that $d\vhi_\alpha=\pi^*dz$.  
The complex structure of $X$ is uniquely determined by this atlas.  
It is easy to check that $\{U_\alpha,g_t^\theta\circ\vhi_\alpha\}$ 
defines a new atlas of charts uniquely determining a new complex 
structure on $X$.  (Here, we used the standard identification 
$\bbC=\bbR^2$.)  The space $X$ with this new complex structure is 
the Riemann surface $X_t^\theta$ referred to in the introduction.  
The family $X_t^\theta$ defines a unit speed geodesic in $\cM_g$ 
with respect to the Teichm\"uller metric $\tau$.  It carries 
a quadratic differential $q_t^\theta$ which is the square the 
holomorphic $1$-form determined by the new charts.  

A \textit{saddle connection} is a geodesic segment which joins a 
pair of branch points without passing through one in its interior.  
Associated to an \textit{oriented} saddle connection $\alpha$ in $X$ 
is a complex number $\int_\alpha\pi^*dz$ which we identify with the 
corresponding vector in $\bbR^2$.  The collection of vectors associated 
to saddle connections in $X$ will be denoted by $V$.  

Let $W=\pm(x_0,y_0)+\bbZ^2$ and $Z=\{(p,q)\in\bbZ^2:\gcd(p,q)=1\}$ 
and note that $$V = W \cup Z.$$  Indeed, a saddle connection in $X$ 
projects to a path in $T$ whose lift to $\bbR^2$ lies in $W\cup Z$.  
Conversely, (\ref{DIO}) implies $\{x_0,y_0,1\}$ is independent over 
$\bbQ$ so that the slope of any vector in $W\cup Z$ is irrational.  
Hence, any vector in $W\cup Z$ can be represented by a geodesic arc 
in $T$ which joins the endpoints of $\gamma$ without passing through 
either one.  The lift of this arc to $X$ is a saddle connection.  

Note that the set of vectors associated with saddle connections 
in $(X_t^\theta,q_t^\theta)$ is simply given by $g_t^\theta V$.  

For any discrete subset $S\subset\bbR^2$, let $\ell(S)$ denote 
the length of the shortest vector in $S$.  To control distances 
in $\cM_g$, we need the following result which is proved in 
slightly greater generality in \cite{Ma93}.  
\begin{Prop}\label{P:M93}
There is a constant $C=C(g)$ such that for all $t\in\bbR$ 
\begin{equation}\label{I:M93}
\tau(X_t^\theta,X_0^\theta)\leq -\log\ell(g_t^\theta V)^2+C.
\end{equation}  
\end{Prop}
Analysing rates of divergence is reduced to studying the function 
$\ell(g_t^\theta V)$, which is carried out in \S\ref{S:Short}.  

\begin{Rem}
The square in (\ref{I:M93}) does not appear in \cite{Ma93} 
due to a different normalisation of the Teichm\"uller metric.  
In our case, the sectional curvature along Teichm\"uller disks 
is $-1$, instead of $-4$.  
\end{Rem}

\subsection{Summable cross products condition}
The surface $X$ carries a flat metric induced by $\pi^*dz$ so that 
it makes sense to talk about parallel lines, area measure, etc.  
For any $\theta\in S^1$, let $F_\theta$ denote the foliation of $X$ 
by lines parallel to $\theta$.  The foliation $F_\theta$ is 
\textit{ergodic} (with respect to area measure) if $X$ cannot be 
written as a disjoint union of two invariant sets of positive measure.  
(An invariant set is one that can be written as a union of leaves.)  
By definition, $\theta$ is a nonergodic direction iff $F_\theta$ 
is not ergodic.  

The next lemma will be useful for finding nonergodic directions.  
\begin{Lem}\label{L:Covers}
Let $\pi':X'\to T$ be the $g$-cyclic branched cover of $T$ along 
another arc $\gamma'$ with the same endpoints as $\gamma$.  Then 
$\pi'$ is biholomorphically equivalent to $\pi$ if and only if 
$\gamma-\gamma'$ represents the trivial element in $H_1(T,\bbZ/g\bbZ)$.  
\end{Lem}
\begin{proof}
Let $U\subset X$ be the set of points lying over the complement 
of $\gamma$ in $T$ and let $U'\subset X'$ be defined similarly.    
We shall identify a dense subset of $U$ with a dense subset of 
$U'$ as follows.  Fix a base point $z_0\not\in\gamma\cup\gamma'$ 
and let $\cU$ be the set of paths in $X$ starting at $z_0$ which 
are transverse to $\pi^{-1}\gamma$.  For any $\alpha\in\cU$, the 
intersection number $i_\gamma(\alpha)\in\bbZ/g\bbZ$ is the number 
of times $\alpha$ crosses $\gamma$ positively.  (This notion depends 
on a choice of orientation for $T$, which we assume has been fixed.) 
The map $\alpha\mapsto(\alpha(1),i_\gamma(\alpha))$ (where $\alpha(1)$ 
denotes the terminal point of $\alpha$) induces a bijection between 
$U$ and $\cU/\sim$ where $$\alpha\sim\alpha'\qquad\text{iff}\qquad 
\alpha(1)=\alpha'(1)\;\;\text{and}\;\;i_\gamma(\alpha)=i_\gamma(\alpha').$$  
Similarly, $U'$ may be identified with classes of paths transverse 
to $\gamma'$ using the above with $\gamma$ replaced by $\gamma'$.  

Note $i_\gamma(\alpha-\alpha')=i_{\gamma'}(\alpha-\alpha')$ iff 
the homology intersection of the cycles $\gamma-\gamma'$ and 
$\alpha-\alpha'$ vanishes.  Therefore, if $\gamma$ is homologous 
to $\gamma'$, there exists a bijection of $U\cap\pi^{-1}\pi'(U')$ 
with $U'\cap\pi'^{-1}\pi(U)$ which extends uniquely to a biholomorphic 
equivalence between $\pi$ and $\pi'$.  Conversely, if $\gamma$ is not 
homologous to $\gamma'$, then there is a closed curve $\beta$ disjoint 
from $\gamma'$ such that $i_\gamma(\beta)\neq0$.  Since its lift is 
closed in $X'$ but not in $X$, $\pi'$ cannot be equivalent to $\pi$.  
\end{proof}

In the sequel, the cross product of two vectors in $\bbR^2$ is 
defined to be a scalar 
$$\vec{a,b}\times\vec{c,d} := ad-bc.$$

\begin{Lem}\label{L:SumX}
Let $(w_j)_{j\geq0}$ be a sequence of vectors of vectors in 
$\pm\vec{x_0,y_0}+g\bbZ^2$ whose Euclidean lengths form an 
increasing sequence and suppose that 
\begin{equation}\label{SumX}
       \sum_{j=0}^\infty |w_j\times w_{j+1}| < \infty.  
\end{equation}
Then $|w_j|^{-1}w_j$ converges to a nonergodic direction in $X$.  
\end{Lem}
\begin{proof}
Note that the direction of the vector $\vec{x_0,y_0}$ associated to 
$\gamma$ is nonergodic because $\pi^{-1}\gamma$ partitions $X$ into 
$g$ invariant sets of equal area.  Similarly, there is associated 
to each $w_j$ a $g$-partition of some branched cover of $T$ that 
is biholomorphically equivalent to $\pi$, by Lemma~\ref{L:Covers}.  
Since a biholomorphism preserves partitions by invariant sets, 
it follows that the direction of each $w_j$ is also nonergodic.  
Now observe that the symmetric difference of the $g$-partitions 
of $X$ associated to a consecutive pair of vectors in the sequence 
is a union parallelograms whose area is bounded above by a constant 
times the absolute cross product.  An elementary argument\footnote
{For more details see the proof of Lemma~1.1 in \cite{Ch}.} shows 
$\lim|w_j|^{-1}w_j$ exists while (\ref{SumX}) implies the sequence 
of $g$-partitions converge measure-theoretically to a $g$-partition 
invariant in the limit direction.  
\end{proof}

\begin{Rem}
Lemma~\ref{L:SumX} is due to Masur-Smillie in genus $2$ and is 
the precursor to a general criterion for nonergodicity developed 
in \cite{MS}.  Their original motivation was to give a geometric 
interpretation, in the context of rational billiards, of certain 
$\bbZ/2$ skew-products studied by Veech in \cite{Ve}.  
\end{Rem}

\section{Analysis of the shortest vector function}\label{S:Short}
The main result of this section is Proposition~\ref{P:PL}, 
which is used to control rates.  Its hypotheses are motivated by 
Lemmas~\ref{L:SVC}, \ref{L:!v} and \ref{L:Nested} while its 
conclusion is motivated by Lemmas~\ref{L:Peak} and \ref{L:Valley}.  

\begin{Lem}\label{L:Peak}
Let $\lambda(t) = -\log |g_t^\theta v|^2$ where $\theta\in S^1$ and 
$v\in\bbR^2$ and assume $|v\times\theta|\neq0$ and $|v\cdot\theta|\neq0$.  
Then the unique maximum $(T,M)$ of $\lambda(t)$ satisfies 
\begin{align}
    T &= \log\frac{|v||v'|}{|v\times v'|} + O(1)\label{I:T}\\
    M &= \log\frac{|v'|}{|v||v\times v'|} + O(1)\label{I:M}
\end{align}
provided $\max(\angle v\theta$, $\angle v'\theta$, $\angle vv')\leq\pi/4$ 
and $|v'\times\theta||v|/|v\times v'|\leq1/2$.  
\end{Lem}
\begin{proof}
From 
\begin{equation}\label{E:Length^2}
|g_t^\theta v|^2 = |v\cdot\theta|^2e^{-t}+|v\times\theta|^2e^t
\end{equation}
we see that $(T,M)$ is given by 
\begin{align}
   T &= \log|v\cdot\theta|-\log|v\times\theta|,           \label{D:T}\\
   M &= -\log|v\cdot\theta| -\log|v\times\theta| - \log 2.\label{D:M} 
\end{align}
Rewrite (\ref{D:T}) by eliminating $|v\cdot\theta|$ in favor of $|v|$ 
to get 
\begin{equation}\label{E:T}
  e^{2T} = \frac{|v\cdot\theta|^2}{|v\times\theta|^2} 
         = \frac{|v|^2}{|v\times\theta|^2}-1.  
\end{equation}
Since $\angle v\theta\leq\pi/4$ iff $T\geq0$, (\ref{E:T}) implies 
$\sqrt{2}|v\times\theta|<|v|$ so that 
\begin{equation}\label{I2:T}
\frac{|v|}{\sqrt{2}|v\times\theta|}\leq e^T \leq \frac{|v|}{|v\times\theta|}.
\end{equation}
Using the triangle inequality, the second hypothesis and the fact that 
the sine function is increasing and nonnegative on $[0,\pi/2]$ we have 
\begin{equation*}
  \frac{|v\times v'|}{|v||v'|} - \frac{|v'\times\theta|}{|v'|} 
             \leq \frac{|v\times\theta|}{|v|} \leq 
  \frac{|v\times v'|}{|v||v'|} + \frac{|v'\times\theta|}{|v'|}. 
\end{equation*}
The first hypothesis now implies 
\begin{equation}\label{I:vv'}
  (1-1/2) \frac{|v\times v'|}{|v||v'|} 
               \leq \frac{|v\times\theta|}{|v|} \leq 
                         (1+1/2) \frac{|v\times v'|}{|v||v'|} 
\end{equation}
which together with (\ref{I2:T}) gives (\ref{I:T}).  
Finally, combining (\ref{D:T}) and (\ref{D:M}) we have 
$e^{-M}=2|v\times\theta|^2e^T$, so that (\ref{I:M}) 
follows from (\ref{I2:T}) and (\ref{I:vv'}).  
\end{proof}

\begin{Lem}\label{L:Valley}
Let $\theta$ and $v$ be as in Lemma~\ref{L:Peak} and assume $v'$ is 
another vector with $|v'\times\theta|\neq0$ and $|v'\cdot\theta|\neq0$.  
If $|v|<|v'|$, $|v'\times\theta|\leq|v\times\theta|/\sqrt{2}$ and 
$|v'\cdot\theta|\geq\sqrt{2}|v\cdot\theta|$ then there exists a 
unique $t>0$ such that $\lambda(t)=-\log|g_t^\theta v'|^2$.  
Moreover, $t$ and $m:=\lambda(t)$ satisfy the following estimates.  
\begin{align}
    t &= \log\frac{|v'|^2}{|v\times v'|} + O(1)\label{I:t}\\
    m &= \log\frac{1}{|v\times v'|} + O(1)\label{I:m}
\end{align}
\end{Lem}
\begin{proof}
From (\ref{E:Length^2}) we see the unique solution to 
$|g_t^\theta v|=|g_t^\theta v'|$ is determined by 
\begin{equation}\label{D:t}
       e^{2t} = \frac{|v'\cdot\theta|^2-|v\cdot\theta|^2}
                     {|v\times\theta|^2-|v'\times\theta|^2}.  
\end{equation}
The second and third hypotheses imply $t$ is well-defined and that 
$e^{2t}>1$ iff $|v|<|v'|$.  Hence, the first hypothesis implies $t>0$.  

Note that the second and third hypotheses hold after $v$ and $v'$ 
are replaced by the vectors $g_t^\theta v$ and $g_t^\theta v'$.  
Since these vectors have the same length, an elementary calculation 
shows the sine of the angle $\phi$ between them is at least $1/3$.  
(Note this is the sine of the angle between $\vec{\sqrt{2},1}$ and 
$\vec{1,\sqrt{2}}$.)  Since $g_t^\theta$ preserves cross products, 
we have $|v\times v'|=e^{-m}\sin\phi$, which implies (\ref{I:m}).  

To get (\ref{I:t}) we first consider the unique maximum time $T'$ 
of the function $t\to-\log|g_t^\theta v'|^2$.  
The analog of (\ref{E:T}) for $T'$ is 
\begin{equation}\label{E:T'}
e^{2T'} = \frac{|v'\cdot\theta|^2}{|v'\times\theta|^2} 
        = \frac{|v'|^2}{|v'\times\theta|^2} - 1.  
\end{equation}
Note that $t\leq T'$ for the second hypothesis together with (\ref{D:t}) 
and (\ref{E:T'}) implies 
$$e^{2t}\leq\frac{|v'\cdot\theta|^2-|v\times\theta|^2}
                 {|v'\times\theta|^2}\leq e^{2T'}.$$  
Now using the definition of $m$, the analog of (\ref{E:Length^2}) for $v'$, 
(\ref{D:t}) and (\ref{E:T'}) we have 
\begin{align*}
e^{-m} &= |v'\cdot\theta|^2e^{-t}+|v'\times\theta|^2e^t \\
       &= |v'|^2e^{-t}\left(\frac{|v'\cdot\theta|^2+|v'\times\theta|^2e^{2t}}
                               {|v'\cdot\theta|^2+|v'\times\theta|^2}\right)\\
       &= |v'|^2e^{-t}\left(1+\frac{e^{2t}-1}{e^{2T'}+1}\right)
\end{align*} 
so that $|v'|^2\leq e^{t-m}\leq2|v'|^2$ since $0<t\leq T$.  
(\ref{I:t}) now follows from (\ref{I:m}).  
\end{proof}

\begin{Lem}\label{L:SVC}
Let $V$ be a discrete subset of $\bbR^2$ such that $\bbR v\cap V=\{\pm v\}$ 
for all $v\in V$ and assume $V\neq\emptyset$.  
Then $\ell(g_t^\theta V)=|g_t^\theta v|$ for some $t\in\bbR$ if 
\begin{equation}\label{I:SVC}
2|v||v\times\theta|\leq
\min\{|v\times u|:|u|\leq\sqrt{2}|v|,u\in V,u\neq\pm v\}.  
\end{equation}
In fact, the condition (\ref{I:SVC}) implies 
$\ell(g_t^\theta V)=|g_t^\theta v|$ for some open interval of $t$ 
near the unique (possibly infinite) time when $|g_t^\theta v|$ is 
minimized.  
\end{Lem}
\begin{proof}
If $v\times\theta=0$ then $|g_t^\theta v|=e^{-t/2}|v|<|g_t^\theta u|$ 
for any $u\in V$ with $|u|>|v|$ since $g_t^\theta$ shrinks Euclidean 
lengths by a factor of at most $e^{t/2}$.  Since $V$ is discrete, 
there are only finitely many $u\in V$ with $|u|\leq|v|$.  For each 
such $u\neq\pm v$ we have $|g_t^\theta u|>|g_t^\theta v|$ for some $t$.  
Therefore, $\ell(g_t^\theta V)=|g_t^\theta v|$ for all large enough $t$.  

If $v\times\theta\neq0$ let $\eps=|g_T^\theta v|$ where $T$ is the 
unique time when $|g_t^\theta v|$ is minimized.  Note that the angle 
between $g_T^\theta v$ and the line $\bbR\theta$ is $\pi/4$ so that 
$$|v\times\theta|e^{T/2}=\eps/\sqrt{2}=|v\cdot\theta|e^{-T/2}.$$  
Suppose $u\in V$ is a vector with $|g_T^\theta u|\leq|g_T^\theta v|$.  
Since $g_T^\theta$ stretches Euclidean lengths by a factor of at most 
$e^{|T|/2}$, we have $|u|\leq\eps e^{|T|/2}\leq\sqrt{2}|v|$ by the 
above equations.  Observing that $g_T^\theta$ preserves cross products, 
we have $$|v\times u|=|g_T^\theta u\times g_T^\theta v|\leq\eps^2
             =2|v\cdot\theta||v\times\theta|<2|v||v\times\theta|$$ 
as $v\times\theta\neq0$.  Now (\ref{I:SVC}) implies $u=\pm v$ so that 
$\ell(g_T^\theta V)=|g_T^\theta v|$.  This proves the first part of 
the lemma while the second part follows by discreteness of $V$.  
\end{proof}

Let $\tW\subset W$ consist of those vectors $w$ for which there is some vector 
$v\in Z$ satisfying $|v|\leq\sqrt{2}|w|$ and $|w\times v|\leq1/2\sqrt{2}$.  
\begin{Lem}\label{L:!v}
For any $w\in\tW$ there exists a unique $v\in Z$ up to sign such that 
$|w\times v|=\min\{|w\times u|:|u|\leq\sqrt{2}|w|,u\in V,u\neq\pm w\}$.  
\end{Lem}
\begin{proof}
Since $W\cap Z=\emptyset$, the hypothesis implies the minimum exists.  
In fact, it must be realized by some vector in $Z$, for if $w'\in W$ 
then $w'=w+du$ for some $u\in Z$ and positive integer $d$ so that 
$|w\times w'|=|w\times du|\geq|w\times u|$.  Now let $v$ be the vector 
associated to $w$ and consider the capped rectangle $R$ defined by the 
inequalities $|u|\leq\sqrt{2}|w|$ and $|w\times u|\leq|w\times v|$.  
It is enough to show $R\cap Z=\{\pm v\}$.  
Suppose there exists $u\in R\cap Z$ such that $u\neq\pm v$.  Note that 
the area of $R$ is $<4\sqrt{2}|w\times v|\leq2$ while the area of the 
parallelogram $P$ with vertices at $\pm u$ and $\pm v$ is exactly $2$.  
This is absurd since $P\subset R$.  Therefore, $R\cap Z=\{\pm v\}$.  
\end{proof}

For any $w\in\tW$ define 
$$I(w)=\{\theta\in S^1: |w\times\theta|<|w\times v|/2|w|,w\cdot\theta>0\}$$  
where $v$ is the vector given by the Lemma~\ref{L:!v}.  

\begin{Lem}\label{L:Nested}
If $w\in\tW$ and $v'\in Z$ satisfy $w\cdot v'>0$ and 
$|w\times v'|\leq1/2\sqrt{2}$, then $w'=w+gv\in\tW$ 
and if $\eps=|w||w\times v'|/|v'||w\times v|\leq1/5$ 
then $\overline{I(w')}\subset I(w)$.  
\end{Lem}
\begin{proof}
Since $|w'\times v'|=|w\times v'|\leq1/2\sqrt{2}$ while $w\cdot v'>0$ 
implies $|w'|>|v'|$, we have $w\in\tW$, easily.  The angle between $w$ 
and $w'$ is at most $2\eps|I(w)|$ (where $|\cdot|$ denotes Lebesgue 
measure induced by arc length) because $|w'|>g|v|$ implies $$\sin\angle 
  ww' = \frac{|w\times w'|}{|w||w'|} < \frac{|w\times v'|}{|w||v'|} 
      = \frac{\eps|w\times v|}{|w|^2} = 2\eps\sin\frac{|I(w)|}{2}$$  
and $\sin^{-1}(2\eps x)\leq2\eps\sin^{-1}x$, $0\leq x\leq1$.  
Similarly, $|I(w')|<g^{-1}\eps|I(w)|$ because 
$$\frac{|w'\times v'|}{|w'|^2} < \frac{|w\times v'|}{g|w'||v'|} 
      = \frac{\eps|w|}{g|w'|} \frac{|w\times v|}{|w|^2}$$  
and $|w|<|w'|$.  Hence, $\overline{I(w')}\subset I(w)$ provided 
$(2+1/2g)\eps\leq1/2$, which follows from $g\geq1$ and $\eps\leq1/5$.  
\end{proof}

\begin{Prop}\label{P:PL}
If $(w_j)_{j\geq0}$ is a sequence in $W$ satisfying 
\begin{enumerate}
\item[(i)] for all $j$, $w_{j+1}=w_j+gv'$ for some $v'\in Z$ 
           with $|v'|>|w_j|$ and $w_j\cdot v'>0$, 
\item[(ii)] $\limsup|w_j\times w_{j+1}|<g/2\sqrt{2}$, and 
\item[(iii)] $\limsup|w_j||w_j\times w_{j+1}|/
                     |w_{j+1}||w_j\times w_{j-1}|<1/(5g+5)$ 
\end{enumerate}
then there exists a piecewise linear function $\Lambda(t)$ satisfying 
$$\limsup_{t\to\infty}|-\log\ell(g_t^\theta V)^2-\Lambda(t)|<\infty$$ 
(with $\theta=\lim|w_j|^{-1}w_j$) and whose critical points are given by 
\begin{align}
  (T_j,M_j) &= \label{E:TjMj}
        \left(\log\frac{|w_j||w_{j+1}|}{|w_j\times w_{j+1}|}, 
              \log\frac{|w_{j+1}|/|w_j|}{|w_j\times w_{j+1}|}\right),\\
  (t_{j+1},m_{j+1}) &= \label{E:tjmj} 
        \left(\log\frac{|w_{j+1}|^2}{|w_j\times w_{j+1}|}, 
              \log\frac{1}{|w_j\times w_{j+1}|}\right).
\end{align}
where $j>j_1$ for some $j_1>0$.  
\end{Prop}
\begin{proof}
First, verify $w_j\in\tW$ for $j$ large enough.  Indeed, by 
(i) $w_{j+1}=w_j+gv'$ for some $v'\in Z$ with $w_j\cdot v'>0$.  
Thus, we easily have $|v'|<\sqrt{2}|w_{j+1}|$ and (ii) implies 
$|w_{j+1}\times v'|=g^{-1}|w_j\times w_{j+1}|\leq1/2\sqrt{2}$.  
Since $v'$ is the vector associated to $w_{j+1}$, we note here 
that for any $\theta\in I(w_{j+1})$ and $j$ large enough 
\begin{equation}\label{I:wtt}
|w_{j+1}\times\theta| < \frac{|w_{j+1}\times v'|}{2|w_{j+1}|} 
     = \frac{|w_j\times w_{j+1}|}{2g|w_{j+1}|}.  
\end{equation}

Next, by Lemma~\ref{L:Nested}, whose hypothesis $\eps\leq1/5$ 
is implied by (iii) and $|w_{j+1}|\leq(g+1)|v'|$ (from (i)), 
we have $\overline{I(w_{j+1})}\subset I(w_j)$ for $j$ large 
enough.  Hence, $\cap_{j\geq j_0}I(w_j)\neq\emptyset$ for some 
$j_0>0$.  Since $|w_{j+1}|>|w_j|$ and $W\subset V$ is discrete, 
we have $\lim|w_j|=\infty$ so that $\lim|I(w_j)|=0$, which implies 
the intersection consists of a single direction and the vectors 
$|w_j|^{-1}w_j\in I(w_j)$ converge to it.  Thus, $\theta$ is 
well-defined; moreover, $\theta\in I(w_j)$ for $j$ large enough.  

To define $\Lambda(t)$, we first note by (iii) there exists a 
$j_1>0$ such that $t_j<T_j$ for all $j>j_1$, while $T_j<t_{j+1}$ 
for all $j\geq0$ since $|w_{j+1}|>|w_j|$.  Let $\Lambda(t)$ be 
the continuous piecewise linear function whose graph is broken 
precisely at the points $(T_j,M_j)$ and $(t_{j+1},m_{j+1})$ 
for $j\geq j_1$;  it is uniquely determined by requiring its 
slope be $+1$ for $t\leq T_{j_1}$.  Hence, each linear piece 
of $\Lambda(t)$ has slope $\pm1$ since $T_j-t_j = M_j - m_j$ 
and $t_{j+1}-T_j = M_j - m_{j+1}$.  

For $j$ large enough we have $\theta\in I(w_{j+1})$ so that 
Lemma~\ref{L:!v} implies (\ref{I:SVC}) holds with $v=w_{j+1}$.  
The hypotheses of Lemma~\ref{L:Peak}, with $w_j$ and $w_{j+1}$ 
in place of $v$ and $v'$, are easily verified using (\ref{I:wtt}) 
and $|w_{j+1}|>|w_j|$.  Hence, we conclude by Lemmas~\ref{L:SVC} 
and \ref{L:Peak} that the points $(T_j,M_j)$ lie within a uniform 
bounded distance of the graph of $f(t):=-\log\ell(g_t^\theta V)^2$.  

Observe that $f(t)$ is $1$-Lipschitz.  Indeed, there exist 
a sequence of vectors $v_k$ in $V$ and a corresponding sequence 
of intervals $I_k$ whose nonoverlapping union is all of $\bbR$ 
such that for all $k$, $f(t)=-\log|g_t^\theta v_k|^2$ for all 
$t\in I_k$.  It is readily seen from (\ref{E:Length^2}) that 
$f'$ is monotone on each $I_k$ with absolute value $\leq1$.  
The proof of the proposition is now complete once we show: 
\textbf{Claim:} the points $(t_j,m_j)$ lie within a 
uniform bounded distance of the graph of $f(t)$.  

To prove the claim we shall apply Lemma~\ref{L:Valley} to 
the vectors $w_j$ and $w_{j+1}$.  The first hypothesis 
$|w_j|<|w_{j+1}|$ follows by (i).  We record here the 
second and third hypotheses for later reference: 
\begin{equation}\label{I:hypo}
|w_{j+1}\times\theta|\leq|w_j\times\theta|/\sqrt{2}\quad\text
{and}\quad|w_{j+1}\cdot\theta|\geq\sqrt{2}|w_j\cdot\theta|.  
\end{equation}
Using (\ref{I:wtt}), the triangle inequality and then 
$|w_{j+1}|>|w_j|$ we have $$|w_{j+1}\times\theta| 
  < \frac{|w_j\times\theta|+|w_{j+1}\times\theta|}{2g}$$ 
which implies the second hypothesis since $g\geq2$.  
Next, using $|w_{j+1}\times\theta|<|w_j\times\theta|$ and 
$|w_{j+1}|\geq\sqrt{2}|w_j|$, which holds by (i) again, 
we obtain the third hypothesis.  
Lemma~\ref{L:Valley} now implies $(t_{j+1},m_{j+1})$ lies 
within bounded distance of the point $(t,m)$ determined by 
   $$e^{-m/2}=|g_t^\theta w_j|=|g_t^\theta w_{j+1}|.$$  
By definition, we have $\ell(g_t^\theta V)\leq e^{-m/2}$.  
To get an inequality in the other direction, let $\phi$ be 
the angle between $g_t^\theta w_j$ and $g_t^\theta w_{j+1}$ 
and $h$ the height of the isosceles triangle formed by them.  
Since $w_{j+1}=w_j+gv'$ we have 
\begin{equation}\label{E:hv'}
h=e^{-m/2}\cos(\phi/2)\qquad\text{and}\qquad 
|g_t^\theta v'|=(2e^{-m/2}/g)\sin(\phi/2).
\end{equation}  
Observe that the distance between any two lines parallel to 
$g_t^\theta v'$ that intersect $g_t^\theta Z$ is an integer 
multiple of $1/|g_t^\theta v'|$ and the same statement holds 
if $W$ is replaced by $Z$.  Hence, the length of any vector in 
$g_t^\theta V$ which is not a multiple of $g_t^\theta v'$ is 
at least $h$, provided $h|g_t^\theta v'|\leq1/2$, but this 
holds because $g_t^\theta$ is area-preserving so that 
$h|g_t^\theta v'|=|w_j\times w_{j+1}|/2g\leq1/2$ by (ii).  
Therefore, $\ell(g_t^\theta V)\geq\min(|g_t^\theta v'|,h)$.  
Now (\ref{I:hypo}) implies $\sin\phi\geq1/3$, i.e. $\phi$ is 
bounded away from $0$ and $\pi$.  Hence, from (\ref{E:hv'}) 
we see there is a universal constant $c>0$ such that 
$\ell(g_t^\theta V)\geq ce^{-m/2}/g$.  It follows that 
$|m-f(t)|\in O(\log g)$ and since $g$ is fixed, this 
completes the proof of the proposition.  
\end{proof}

\section{Density of primitive lattice points}\label{S:Dens}
The main result of this section is Corollary~\ref{C:Dens2}.  
It will be needed in \S\ref{S:NEsub} to find, given a vector 
$w\in W$, vectors $w'\in W$ such that $w'=w+gv$ for some $v\in Z$ 
satisfying certain given inequalities on $|v|$ and $|w\times v|$.  

\subsection{Continued fractions in vector form}
Recall each $\alpha\in\bbR$ admits an expansion of the form 
\begin{equation}\label{Cfrac}
\alpha = a_0 + \cfrac{1}{a_1+\cfrac{1}{a_2+\dots}}
\qquad a_0\in\bbZ, a_1,a_2,\dots\in\bbN
\end{equation}
whose terms are uniquely determined except for a two-fold 
ambiguity when $\alpha$ is rational; e.g. $22/7=3+1/7=3+1/(6+1/1)$.  
The $k$th \textit{convergent} of $\alpha$ is the reduced fraction 
$p_k/q_k$ that results upon simplifying the expression obtained by 
truncating (\ref{Cfrac}) so that the last term is $a_k$.  
The convergents of $\alpha$ satisfy the recurrence relations 
\begin{equation}\label{RR}
\begin{tabular}{cc}
$p_{k+1} = a_{k+1}p_k + p_{k-1}$ & $\;p_0=a_0,\,p_{-1}=1$\\
$q_{k+1} = a_{k+1}q_k + q_{k-1}$ & $\;q_0=1,\;\;q_{-1}=0$
\end{tabular},
\end{equation}
the identity 
\begin{equation}\label{PR}
      p_kq_{k+1}-p_{k+1}q_k=(-1)^{k+1}
\end{equation}
and the inequalities 
\begin{equation}\label{CFI}
\frac{1}{q_k(q_{k+1}+q_k)} 
        < \left|\alpha-\frac{p_k}{q_k}\right|
        < \frac{1}{q_kq_{k+1}}.  
\end{equation}
A rational $p/q$ is said to be a \textit{best approximation of the 
second kind} if 
\begin{equation}\label{BA2}
|q\alpha-p|\leq|n\alpha-m|\qquad\text{for all}\;\;m\in\bbZ,\;n=1,\dots,q
\end{equation}
and this property characterises the convergents of $\alpha$ modulo 
the $0$th convergent $a_0$, which is a best approximation to $\alpha$ 
iff the fractional part of $\alpha$ is $\leq1/2$.  The following is 
a useful test for a rational to be a convergent: 
\begin{equation}\label{SCC}
\left|\alpha-\frac{p}{q}\right| \leq \frac{1}{2q^2},\;\;\gcd(p,q)=1 
\quad\Rightarrow\quad (p,q)=\pm (p_k,q_k) \;\text{for some}\;k\geq0.  
\end{equation}

It will be convenient for us to recast the above facts in vector form.  
Setting $v_k:=\vec{p_k,q_k}$, the recurrence relations (\ref{RR}) and 
the identity (\ref{PR}) become 
\begin{gather}\label{VRR}
   v_{k+1} = a_{k+1}v_k + v_{k-1}\qquad 
               v_0=\vec{a_0,1}, \quad v_{-1}=\vec{1,0} \\
\intertext{and}\label{VPR}
   v_k\times v_{k+1} := p_kq_{k+1}-p_{k+1}q_k = (-1)^{k+1}.  
\end{gather}
Although (\ref{CFI}) can easily be rewritten in vector notation, the 
resulting expression looks awkward because of the distinguished nature of 
the coordinate directions.  Instead, we shall use the following analog of 
(\ref{CFI}) which is expressed in terms of the vector $w:=\vec{\alpha,1}$ 
and its Euclidean length $|w|$:  
\begin{equation}\label{VCFI}
\frac{1}{|v_{k+1}+v_k|} < \frac{\left|w\times v_k\right|}{|w|} 
                        < \frac{1}{|v_{k+1}|}  
\end{equation}
To see (\ref{VCFI}) recall that convergents alternate on both 
sides of $\alpha$ and (\ref{CFI}) follows from the fact 
that the rational $(p_k+p_{k+1})/(q_k+q_{k+1})$ always lies on the 
same side of $\alpha$ occupied by $p_k/q_k$.  (\ref{VCFI}) follows 
similarly from a comparison of the components of $v_{k+1}+v_k$, $w$ 
and $v_{k+1}$ in the direction perpendicular to $v_k$.  

\begin{Def}\label{D:Spec}
If $\theta=|w|^{-1}w$ where $w=\vec{\alpha,1}$ as above, then 
we define $$\Spec(\theta)=\{v_0,v_1,\dots\} = \{v_k\}_{k\geq0}$$  
and call the vectors $v_k$ the \textit{convergents} of $\theta$.  
The definition is extended to all unit vectors by requiring 
$$\Spec(1,0)=\{\vec{1,0}\},\quad \Spec(-\theta)=-\Spec(\theta)$$ 
and to all nonzero vectors by $\Spec(w)=\Spec(|w|^{-1}w)$.  
We shall also denote by $\spec(w)$ the sequence of Euclidean 
lengths of vectors in $\Spec(w)$.  
\end{Def}

The next lemma was motivated by (\ref{SCC}) and will be needed 
in \S\ref{S:NEsub}.  
\begin{Lem}\label{L:VSCC}
Let $w$ be a vector that makes an angle $\phi$ with the $y$-axis.  
Then for any $v\in\bbZ^2$ such that $|v|\cos\phi>1$ we have 
\begin{equation}\label{VSCC}
\frac{|w\times v|}{|w|} \leq \frac{1}{2|v|},\quad\gcd(v)=1. 
\quad\Rightarrow\quad \pm v\in\Spec(w).  
\end{equation}
\end{Lem}
\begin{proof}
Let $P=P(v,w)$ be the closed parallelogram that has $\pm v$ as two 
of its vertices, one pair of sides parallel to $w$ and the other pair 
parallel to the $x$-axis.  The characterisation of convergents given 
in (\ref{BA2}) is equivalent to the statement that every nonzero 
$u\in P\cap\bbZ^2$ belongs to the union of the two sides of $P$ 
parallel to $w$.  Hence, it is enough to verify this statement under 
the given hypotheses.  

Apply Lemma~\ref{L:SVC} with $\theta=|w|^{-1}w$ (and $V$ the set 
of integer lattice points which are not scalar multiples of $v$) 
to conclude there is some $T$ for which $g_T^\theta v$ is the 
shortest vector in $V'=g_T^\theta(\bbZ^2-0)$.  Let $E$ be the 
inverse image under $g_T^\theta$ of the largest closed disk 
centered at the origin whose interior is disjoint from $V'$.  
The boundary of $E$ is an ellipse passing through the points 
$\pm v$ while the interior contains no integer lattice points 
other than the origin.  

Without loss of generality, we assume $w$ lies in the first 
quadrant and $v$ in the upper half plane.  There are two cases.  
First, if $v$ lies to the right of $w$, then $E$ contains $P$ 
and we are done.  Now, if $v$ lies to the left of $w$, 
then let $x$ be the length of a horizontal side of $P$ and 
$y$ the vertical distance between $v$ and its reflection in 
the line $\bbR w$.  It is enough to show $x<1$ and $y<1$.  
If $z$ is the distance between $v$ and its reflection then 
$z=2|v|\sin\angle{wv}=2|w\times v|/|w|\leq1/|v|<\cos\phi$ 
so that $x=z\sec\phi<1$ and $y=z\sin\phi<1$.  
\end{proof}

\subsection{Density of rationals in intervals}
For any $\Omega\subset\bbR^2$ with $0<\area(\Omega)<\infty$ let 
$$\dens(\Omega) = \frac{^\#Z\cap\Omega}{\area(\Omega)}.$$  
\begin{Lem}\label{L:8/27}
If $\Omega$ is a compact convex subset of the first quadrant 
containing $(0,0)$, $(1,0)$ and $(0,1)$ but not $(1,1)$ then 
$\dens(2\Omega\setminus\Omega)>8/27$.  
\end{Lem}
\begin{proof}
By convexity there is a function $y=f(x), 0\leq x\leq1$ whose graph 
is contained in $\partial\Omega$ and $f(0)\geq1$.  Similarly, there is 
a $x=g(y), 0\leq y\leq1$ whose graph is contained in $\partial\Omega$ 
and $g(0)\geq1$.  Without loss of generally, assume $$f(1/2)\geq g(1/2).$$  
Let $\Omega_1\subset\Omega$ be the part below the graph of $f$ and 
$\Omega_2=\Omega\setminus\Omega_1$.  Since $\Omega_1$ lies below 
any tangent line at $(1/2,f(1/2))$, $\area(\Omega_1)\leq f(1/2)$.  
There are three cases.  First, if $f(1/2)<1$ then the assumption 
above implies $g(0)<3/2$ so that $\area(\Omega_2)\leq1/8$.  
Since $(1,1)\in2\Omega\setminus\Omega$ we have 
$$\dens(2\Omega\setminus\Omega)\geq 
\frac{1}{3}\left(\frac{1}{f(1/2)+1/8}\right)>8/27.$$
Second, if $f(1/2)\geq1$ and $f(1)<1/2$ then it still follows that 
$g(0)<3/2$ while the number of vectors in $2\Omega\setminus\Omega$ 
of the form $(1,n)$ is $\lfloor2f(1/2)\rfloor$ so that 
$$\dens(2\Omega\setminus\Omega)>
\frac{1}{3}\left(\frac{2f(1/2)-1}{f(1/2)+1/8}\right)\geq8/27.$$
Finally, if $f(1/2)\geq1$ and $f(1)\geq1/2$ then the assumption above 
implies $g(0)<2$ so that $\area(\Omega_2)\leq1/2$.  Apart from points of 
the form $(1,n)$ we also have $(2,1)\in2\Omega\setminus\Omega$.  Hence, 
$$\dens(2\Omega\setminus\Omega)>
\frac{1}{3}\left(\frac{2f(1/2)}{f(1/2)+1/2}\right)\geq4/9>8/27.$$
\end{proof}

Note if $\Omega_a=\{(x,y): x+y\leq a, x\geq0, y\geq0\}$ then 
\begin{align*}
\dens(2\Omega_a\setminus\Omega_a) &= 2/3a^2 
             \quad\text{for}\quad 1\leq a<3/2\quad\text{and}\\
\dens(\gamma\Omega_1\setminus\Omega_1) &= 0
             \quad\text{for}\quad \gamma<2
\end{align*}
show that the constants in the preceding lemma are sharp.  

Let $S^1(\bbQ)$ be the set of unit vectors of the form $v/|v|$ for some 
$v\in Z$ and $S^1_b(\bbQ)$ the subset formed by those with $|v|\leq b$.  
Let $\cI_b$ denote the collection of intervals in $S^1$ with endpoints 
in $S^1_b(\bbQ)$.  For any interval $I\subset S^1$ let 
$$\Omega(I,b) = \{v\in\bbR^2 : \bbR_+v\cap I\not=\emptyset, |v|\leq b\}.$$

\begin{Prop}\label{P:8/27}
For any $I\in\cI_b$, $\dens(2\Omega(I,b)\setminus\Omega(I,b))>8/27$.  
\end{Prop}
\begin{proof}
First we show if $I$ is minimal, i.e. $S^1_b(\bbQ)\cap\text{int}I=\emptyset$, 
then its endpoints correspond to a pair of vectors $v,v'\in Z$ such that 
$|v\times v'|=1$.  Indeed, there is a linear map $\gamma$ that sends $v$ 
to $(0,1)$ and $v'$ to $(a,a')\in\bbZ^2$ with $0\leq a'<a=|v\times v'|$.  
If $a>1$ then $\gamma^{-1}(1,1)\in S^1_b(\bbQ)\cap\text{int}I$; hence, $a=1$.  
Now observe $\gamma\Omega(I,b)$ is a compact convex set satisfying the 
hypothesis of Lemma~\ref{L:8/27} and since $\gamma$ preserves density, 
the proposition holds for minimal $I$ in $\cI_b$.  Since every interval 
in $\cI_b$ is a (finite) disjoint union of minimal ones, this completes 
the proof.  
\end{proof}


\begin{Thm}\label{T:Dens}
  Let $\Omega=\Omega(I,b)$ where 
  $I = \{\theta'\in S^1 : \sin\angle\theta\theta'<\eps/b\}$, 
  $\theta\in S^1$, $\eps>0$ and $b\geq1$.  
  Then $\spec(\theta)\cap[\eps^{-1},b]\not=\emptyset$ 
  implies $\dens(2\Omega\setminus\Omega)>4/27\pi$.  
\end{Thm}
\begin{proof}
Let $v_k\in\Spec(\theta)$ be the convergent with length 
$|v_k|=\max\spec(\theta)\cap[\eps^{-1},b]$.  Then the RHS of 
(\ref{VCFI}) implies (the direction of) $v_k$ lies in $I$.  
Without loss of generality, assume $v_k$ lies to the left of $\theta$.  
Let $\theta'$ be the right endpoint of $I$ and $v'_l\in\Spec(\theta')$ 
the convergent with length $|v'_l|=\max\spec(\theta)\cap[1,b]$.  

Note that $v_k$ does not lie strictly between $v'_l$ and $v'_{l+1}$ 
since the length of any such vector in $Z$ is at least $|v'_{l+1}|>b$.  
Nor can $v_k=v'_l$ since the RHS of (\ref{VCFI}) would imply $v'_l$ 
lies strictly to the right of $v_k$.  Therefore, $v'_l$ lies strictly 
to the right of $v_k$.  Let $\Omega'=\Omega(J,b)$ where $J$ is 
the interval with left endpoint $v_k$ and right endpoint $v'_l$.  

The interval $I$ has length $|I|=2\sin^{-1}(\eps/b)\leq\pi\eps/b$.  
If $|v'_l|\geq2\eps^{-1}$ then the RHS of (\ref{VCFI}) implies 
$|J| > \sin^{-1}(\eps/b) - \sin^{-1}(\eps/2b) > \eps/2b$ (we may 
assume $\eps<b$ for otherwise the theorem is easily seen to hold) 
while if $|v'_l|<2\eps^{-1}$ then 
$|J| > \sin\angle v_kv'_l \geq \frac{1}{|v_k|\,|v'_l|} > \eps/2b.$  
In either case, we have $|J|>|I|/2\pi$.  

If $v'_l$ lies to the left of $\theta'$ or $|v'_{l+1}|>2b$ then 
$Z\cap(2\Omega'\setminus\Omega')\subset Z\cap(2\Omega\setminus\Omega)$ 
since any vector strictly between $v'_l$ and $v'_{l+1}$ has length 
greater than $2b$.  In this case, Proposition~\ref{P:8/27} implies
$\dens(2\Omega\setminus\Omega)>4/27\pi$.  If $|v'_l|\geq2\eps^{-1}$ 
then arguing as before we see the angle between $v_k$ and $v'_{l+1}$ 
is at least $\eps/2b>|I|/2\pi$ so that we may again conclude that 
$\dens(2\Omega\setminus\Omega)>4/27\pi$.  Therefore, we may assume 
$v'_{l+1}$ lies between $v_k$ and $\theta'$, $|v'_{l+1}|\leq2b$ and 
$|v'_l|<2\eps^{-1}$ in the remaining.  Assuming $\eps<b/\sqrt{2}$ 
as we may, since the theorem is easily seen to hold otherwise, 
we obtain the following criterion for a vector left of $\theta'$ 
to lie in $I$: 
\begin{equation}\label{I:vinI}
\frac{|v\times\theta'|}{|v|} < \frac{\sqrt{2}\eps}{b}
\end{equation}
Note that a vector of the form $av'_l+v_{l-1}$ lies to the left of 
$\theta'$ if $a\leq a_{k+1}$.  We will show the number of vectors 
of length $\geq b/2$ satisfying the above conditions is at least 
$c_0b\eps$ for some absolute constant $c_0>0$.  This will complete 
the proof since between any two vectors of length $\geq b/2$ (but 
$\leq b$) there is a vector whose length is $>b$ and $\leq2b$.  
Using the RHS of (\ref{VCFI}) we have 
\begin{equation}\label{I:Ifrac}
|v\times\theta'| = |v'_{l+1}\times\theta'-(a_{l+1}-a)v'_l\times\theta'|
               \leq \frac{1+|a_{l+1}-a|}{|v'_{l+1}|}.  
\end{equation}
Assuming $|v|\geq b/2$ ($v=av'_l+v_{l-1}$) we see (\ref{I:vinI}) 
holds for $\lfloor b\eps/\sqrt{2}\rfloor$ integers $a\leq a_{l+1}$.  
Among these there are at least $b/2|v'_l|\geq b\eps/4$ which make 
$|v|\geq b/2$.  Since $\area(\Omega)<\pi b\eps/2$, we conclude 
$\dens(2\Omega\setminus\Omega)>1/6\pi$.  
\end{proof}

\subsection{Density in a strip}
\begin{Thm}\label{T:Dens2}
  Let $\Sigma=\{v\in\bbR^2:|v\times\theta|<\eps,b<|v|\leq2b,v\cdot\theta>0\}$ 
  where $\theta\in S^1$, $0<\eps\leq1$ and $b\geq1$.  
  Then $\spec(\theta)\cap[\eps^{-1},b]\neq\emptyset$ 
  implies $\dens(\Sigma)>2/27\pi$.  
\end{Thm}
\begin{proof}
Let $v_k\in\Spec(\theta)$ with $|v_k|=\max\spec(\theta)\cap[\eps^{-1},b]$.  
There are two cases.  
If $|v_k|\geq2\eps^{-1}$ let $\Omega=\Omega(\theta,\eps,b)$ be the 
region in Theorem~\ref{T:Dens} and put $\Omega'=\Omega(\theta,\eps/2,b)$.  
Then Theorem~\ref{T:Dens} implies $\dens(2\Omega'\setminus\Omega')>4/27\pi$.  
Since $2\Omega'\setminus\Omega'\subset\Sigma$ and $\area(\Omega')>b\eps/2$ 
it follows there are at least $2b\eps/9\pi$ vectors in $Z\cap\Sigma$.  
If $|v_k|<2\eps^{-1}$ let $v=av_k+v_{k-1}$ and use (\ref{I:Ifrac}) for 
$\theta$ and the hypothesis $|v_k|\geq\eps^{-1}$ to deduce 
$$|v\times\theta| \leq \frac{1+|a_{k+1}-a|}{|v_{k+1}|} 
                  < \frac{1+|a_{k+1}-a|}{a_{k+1}|v_k|} \leq \eps$$
for all $a$ such that $1\leq a<2a_{k+1}$.  Since $|v_k|<2\eps^{-1}$ 
the number $v$ with $b<|v|\leq2b$ is $\geq b\eps/2$.  
In either case, we have found $2b\eps/9\pi$ vectors $Z\cap\Sigma$.  
By an elementary analysis we obtain 
\begin{equation}\label{I:area}
3b\eps/2\leq 2b\eps\left(1-\frac{\eps^2}{4b^2}\right) \leq 
\area(\Sigma) \leq 2b\eps\left(1+\frac{\eps^2}{2b^2}\right) \leq 3b\eps.  
\end{equation}
Using the RHS we get $\dens(\Sigma)>2/27\pi$.  
\end{proof}

\begin{Cor}\label{C:Dens2}
There are universal constants $\rho_1>0$ and $c_1>0$ such that 
$\spec(w)\cap [\eps^{-1},b]\neq\emptyset$ (and $0<\eps\leq1\leq b$) 
implies there are $\rho_1b\eps$ vectors $v\in Z$ satisfying the 
inequalities $$w\cdot v>0,\quad 
b\leq|v|\leq2b,\quad c_1\eps|w|\leq|w\times v|\leq\eps|w|.$$  
\end{Cor}
\begin{proof}
Let $\Sigma=\Sigma(\eps,b)$ be as in Theorem~\ref{T:Dens2} with 
$\theta=|w|^{-1}w$.  First we claim $\dens(\Sigma)$ is bounded 
above by some universal constant.  Indeed, from the cross product 
formula we see there's a universal constant $C>1$ such that the 
largest (resp. smallest) angle between two vectors in $Z\cap\Sigma$ 
is $<C\eps/b$ (resp. $>C^{-1}/b^2$).  Hence, the number of vectors 
in $Z\cap\Sigma$ is at most $C^2b\eps$ so that the LHS of (\ref{I:area}) 
implies the claim.  

Now let $\Sigma_1=\Sigma(c_1\eps,b)$ and $\Sigma_2=\Sigma\setminus\Sigma_1$.  
Observe that the claim implies for any $\rho<2/27\pi$, $c_1$ can be chosen 
sufficiently small so that $\dens(\Sigma_2)>\rho$.  Since (\ref{I:area}) 
implies $\area(\Sigma_2)\geq cb\eps$ for some universal $c>0$, the corollary 
follows.  
\end{proof}

\section{Nonergodic directions and sublinear growth}\label{S:NEsub}
In this section we prove Theorem~1.  

Let $e_0>\max(d_0,2)$ be given.  
The construction involves the choice of a sequence $(\delta_j)_{j\geq0}$ 
descending to zero at some prescribed rate as required by 
Lemmas~\ref{L:Sublinear} and \ref{L:eps} below.  
For concreteness we set 
\begin{equation*}
    \delta_j:=\frac{e_0}{j+1}\quad\text{for all $j\geq0$.}  
\end{equation*}
In addition, we also fix a constant 
\begin{equation*}
    C:=\max(2g+1,c_1^{-1}e^{2e_0})
\end{equation*}
needed in the statement of Lemma~\ref{L:Good} below.  

Our goal is to find sequences $(w_j)_{j\geq0}$ in $W$ satisfying 
\begin{equation}\label{I:Child1}
  |w_j|^{1+\delta_j}\leq |w_{j+1}|\leq C|w_j|^{1+\delta_j}, \quad
  \frac{C^{-1}}{\log|w_j|}\leq|w_j\times w_{j+1}|\leq\frac{C}{\log|w_j|}
\end{equation}
and 
\begin{equation}\label{I:Child2}
  w_{j+1}=w_j+gv'\quad\text{for some $v'\in Z$ with}\quad
  |v'|>|w_j|\;\text{and}\; w_j\cdot v'>0
\end{equation}
for all $j\geq0$.  
\begin{Def}\label{D:Child}
If $w_j$ and $w_{j+1}$ satisfy (\ref{I:Child1}) and (\ref{I:Child2}) 
then we say $w_{j+1}$ is a \textit{child} of $w_j$.  More precisely, 
(\ref{I:Child1}) and (\ref{I:Child2}) define a family $\{\prec_j\}$ 
of binary relations on $W$ and to say $w_{j+1}$ is a child of $w_j$ 
is equivalent to the statement $w_j\prec_j w_{j+1}$.  
\end{Def}

\begin{Def}\label{D:Admissible2}
We say $(w_j)_{j\geq0}$ is \textit{admissible} if $|w_0|>1$, 
$w_0\pm(x_0,y_0)\in g\bbZ^2$ and for all $j\geq0$, $w_{j+1}$ 
is a child of $w_j$.  A finite sequence $(w_0,\dots,w_k)$ is 
admissible if the latter condition holds for $0\leq j<k$.  
\end{Def}

The choice of $(\delta_j)$ was motivated by the next two lemmas, 
which are stated more generally for a sequence of positive $\delta_j$.  
\begin{Lem}\label{L:Sublinear}
  If $(\delta_j)_{j\geq0}$ is a sequence of positive real numbers such that 
$$\liminf j\delta_j>1\quad\text{and}\quad e_0=\limsup j\delta_j<\infty$$ 
  and $(w_j)_{j\geq0}$ is an admissible sequence then 
  $\lim|w_j|^{-1}w_j$ is a slowly divergent nonergodic direction 
  whose sublinear rate is at most $1-1/e_0$.  
\end{Lem}
\begin{proof}
First note that the hypotheses imply 
(i) $\lim\delta_j=0$ and (ii) $\sum\delta_j=\infty$.  
Let $S_j=\sum_{i<j}\delta_i$ and $R_j=\sum_{i<j}\log(1+\delta_i)$.  
Claim: $\lim R_j/S_j=1$.  
Indeed, for any $x\in[0,1]$ we have $x(1-x)\leq\log(1+x)\leq x$.  
Using (i) we may fix $K>1$ so that $R_j\leq S_j\leq K R_j$ for all $j$.  
Then $\lim R_j=\infty$ by (ii).  
Let $c>1$ be given.  Using (i) again we fix $j_0$ large enough 
so that $\delta_j\leq\sqrt{c}\log(1+\delta_j)$ for all $j\geq j_0$.  
Now $S_j\leq K R_{j_0} + \sqrt{c}R_j\leq c R_j$ for all large 
enough $j$.  Since $c>1$ was arbitrary, this poves the claim.  

From the first inequality in (\ref{I:Child1}) we have 
$\log|w_j|\geq(\log|w_0|)\prod_{i<j}(1+\delta_i)$.  
Now $\liminf j\delta_j>1$ implies for some $p>1$ and $C'>0$ we have 
$S_j\geq p\log j - C'$ for all $j>0$.  The same statement for $R_j$ 
holds by the preceding claim.  
Hence, $\prod_{i<j}(1+\delta_i)\geq c_1j^p$ for some $c_1>0$.  
The last two inequalities in (\ref{I:Child1}) imply the cross 
products form a summable series while (\ref{I:Child2}) implies 
$|w_j|$ is increasing.  Hence, $\lim|w_j|^{-1}w_j$ is a nonergodic 
direction, by Lemma~\ref{L:SumX}.  

It is clear from the preceding that $\lim|w_j|=\infty$; 
however, a much stronger statement holds.  First, there 
is some $p>1$ and $c_2>0$ such that (for all $j>0$) 
\begin{equation}\label{I:|w_j|>}
  \delta_j\log|w_j| \geq (c_2\log|w_0|)j^{p-1}.
\end{equation}
Now the second inequality in (\ref{I:Child1}) implies 
$\log|w_j|\leq(\log|w_0|)\prod_{i<j}(1+\delta_i)+j\log C$.  
Using $\limsup j\delta_j<\infty$ and arguing as before we 
find $q>e_0$ and $C''>0$ such that $R_j\leq q\log j + C''$ 
for all $j>0$.  Thus $\prod_{i<j}(1+\delta_i)\leq c_3j^q$ 
some $c_3>1$ and $\log|w_j|\leq c_3j^q\log|w_0|+j\log C$.  
Using $\log(1+x+y)\leq\log(1+x)+\log(1+y)$ we conclude: 
for some $q\geq p$ and $C'''>0$ we have (for all $j>0$) 
\begin{equation}\label{I:|w_j|<}
\log\log|w_j| \leq \log\log|w_0| + (q+1)\log j + C'''.  
\end{equation}
It follows from (\ref{I:|w_j|>}) and (\ref{I:|w_j|<}) that 
$\lim|w_j|^{\delta_j}/\log|w_j|=\infty$.  

The hypotheses of Proposition~\ref{P:PL} are satisfied since 
(i) is the same as (\ref{I:Child2}) while (\ref{I:Child1}) 
and (\ref{I:|w_j|>}) imply $\lim |w_j|/|w_{j+1}|=0$ and 
$\lim|w_j\times w_{j+1}|=0$, where the ratio of consecutive 
cross-product terms is bounded.  Propositions~\ref{P:M93} and 
\ref{P:PL} imply $\lim|w_j|^{-1}w_j$ is slowly divergent since 
$$\lim \frac{M_j}{T_j} 
       = \lim \frac{\delta_j\log|w_j|-\log|w_j\times w_{j+1}|}
                   {(2+\delta_j)\log|w_j|-\log|w_j\times w_{j+1}|}
       = \lim \frac{\delta_j}{2+\delta_j}=0$$ 
while the sublinear rate is at most 
\begin{align*}
\limsup \frac{\log M_j}{\log T_j} 
       &= \limsup \frac{\log(\delta_j\log|w_j|-\log|w_j\times w_{j+1}|)}
                   {\log((2+\delta_j)\log|w_j|-\log|w_j\times w_{j+1}|)}\\
       &= 1-\liminf \frac{-\log\delta_j}{\log\log|w_j|} 
       \leq 1 - \frac{1}{q}
\end{align*}
because $-\log\delta_j\geq\log j+O(1)$ and $R_j\leq q\log j+O(1)$.  
The proof is completed by observing that $q>e_0$ may be chosen 
arbitrarily close to $e_0$. 
\end{proof}

\begin{Lem}\label{L:eps}
  Let $(\delta_j)_{j\geq0}$ be a sequence of positive real numbers such that 
  $$\liminf j\delta_j>2\quad\text{and}\quad \limsup j\delta_j<\infty$$ 
  and $(w_j)_{j\geq0}$ an admissible sequence.  Then 
  for any $\eps>0$ there exists $L_0=L_0(\eps)$ such that 
  $|w_0|\geq L_0$ implies 
\begin{equation}\label{I:eps}
  \sup_{j\geq0}\frac{(\log |w_j|)^2}{|w_j|^{\delta_j\delta_{j+1}}}\leq\eps.  
\end{equation}
\end{Lem}
\begin{proof}
Repeating the arguments in the preceding proof with the stronger 
hypotheses we find there are constants $p>2$, $c_2>0$, $q\geq p$ 
and $C'''>0$ such that (\ref{I:|w_j|>}) and (\ref{I:|w_j|<}) hold 
for all $j>0$.  It follows that the difference 
\begin{align*}
\delta_j\delta_{j+1}\log|w_j|-2\log\log|w_j| &\geq \\
    (c_2'\log|w_0|)j^{p-2}&-2\log\log|w_0|-2(q+1)\log j-2C'''
    =: \beta(j)
\end{align*}
for some $c_2'>0$.  The function $\beta(j)$ is increasing for $j\geq1$ 
and by choosing $L_0$ large enough we have $\beta(1)\geq-\log\eps$.  
By choosing $L_0$ even larger so that $(\log|w_0|)^2/|w_0|^{\delta_0
\delta_1}\leq\eps$ we obtain (\ref{I:eps}).  
\end{proof}

\begin{Lem}\label{L:Good}
Let $(w_0,\dots,w_j)$ be a admissible sequence and suppose 
\begin{equation}\label{D:Good}
\spec(w_j)\cap[e^t|w_j|\log|w_j|,|w_j|^{1+\delta_j}]\neq\emptyset 
\end{equation}
for some $t\in[0,2e_0]$.  
Then $w_j$ has at least $\rho_1e^{-2e_0}|w_j|^{\delta_j}/\log|w_j|$ 
children  
and these vectors satisfy 
\begin{equation}\label{NotSoGood}
\spec(w_{j+1})\cap[e^{t-\delta_j}
|w_{j+1}|\log|w_{j+1}|,|w_{j+1}|^{1+\delta_{j+1}}]\neq\emptyset
\end{equation}  
provided $|w_0|$ is large enough.
\end{Lem}
\begin{proof}
Apply Corollary~\ref{C:Dens2} with $\eps^{-1}=e^t|w_j|\log|w_j|$ and 
$b=|w_j|^{1+\delta_j}$ to get $\rho_1e^{-t}|w_j|^{\delta_j}/\log|w_j|$ 
vectors $v\in Z$ satisfying the inequalities 
\begin{equation}\label{I:v}
w_j\cdot v>0,\quad|w_j|^{1+\delta_j}\leq|v|\leq 2|w_j|^{1+\delta_j},\quad
\frac{c_1e^{-t}}{\log|w_j|}\leq|w_j\times v|\leq\frac{e^{-t}}{\log|w_j|}.
\end{equation}
The vector $w_{j+1}=w_j+gv$ satisfies (\ref{I:Child2}) by the first 
inequality in (\ref{I:v}), which together with the second inequality 
implies $|w_{j+1}|\geq|v|\geq|w_j|^{1+\delta_j}$.  The third implies 
$|w_{j+1}|\leq|w_j|+g|v|\leq(2g+1)|w_j|^{1+\delta_j}$, which together 
with the remaining inequalities and $|w_j\times w_{j+1}|=g|w_j\times v|$ 
implies (\ref{I:Child1}) for the given value of $C$.  
Therefore, $w_{j+1}$ is a child of $w_j$ and since $t\leq2e_0$, 
this proves the first part.  

Using $|w_{j+1}|>g|v|$, the last inequality in (\ref{I:v}) and 
$g\geq2$ we have  
$$\frac{|w_{j+1}\times v|}{|w_{j+1}|} < \frac{|w_j\times v|}{g|v|} 
       \leq \frac{e^{-t}}{g|v|\log|w_0|} < \frac{1}{2|v|}$$ 
as soon as $\log |w_0|>1$.  Since $v\in Z$, Lemma~\ref{L:VSCC} 
implies $+v$ is a convergent of $w_{j+1}$, where the sign follows 
from $w_{j+1}\cdot v>0$ and the fact that all angle between a vector 
and its convergents are acute.  (See Remark~\ref{R:VSCC} below for 
an explanation of how the hypothesis of Lemma~\ref{L:VSCC} is satisfied.)  

Let $v'\in\Spec(w_{j+1})$ be the next convergent after $v$.  
Using the RHS of (\ref{VCFI}), the second to last inequality 
in (\ref{I:v}), $t\leq2e_0$, $|w_j|<|w_{j+1}|$ and (\ref{I:eps}) 
we have $$|v'|\leq|w_{j+1}||w_j\times v|^{-1}\leq 
c_1^{-1}e^{2e_0}|w_{j+1}|\log|w_{j+1}|\leq|w_{j+1}|^{1+\delta_{j+1}}$$ 
provided $|w_0|$ is chosen large enough as required by 
Lemma~\ref{L:eps} for $\eps=c_1e^{-2e_0}$.  

Using the LHS of (\ref{VCFI}), $|v|<|w_{j+1}|$, $t\geq0$ and 
$\delta_j\leq\delta_0$ we have 
\begin{align*}
|v'|&\geq|w_{j+1}||w_j\times v|^{-1}-|v|\geq|w_{j+1}|(e^t\log|w_j|-1)\\
&\geq e^{t-\delta_j}|w_{j+1}|(e^{\delta_j}\log|w_j|-e^{\delta_0}).  
\end{align*} 
On the other hand, since 
$|w_{j+1}|\leq|w_j|+g|v|\leq(2g+1)|w_j|^{1+\delta_j}$ we have 
\begin{align*}
\log|w_{j+1}| 
              &\leq e^{\delta_j}\log|w_j|-\delta_j^2\log|w_j|+\log(2g+1)
\end{align*}
from which is follows $|v'|\geq e^{t-\delta_j}|w_{j+1}|\log|w_{j+1}|$ 
if $|w_0|$ is chosen large enough as required by Lemma~\ref{L:eps} 
for $\eps= (2g+1)^{-1}e^{-e^{\delta_0}}$.  
\end{proof}

\begin{Rem}\label{R:VSCC}
In order to satisfy the hypothesis of Lemma~\ref{L:VSCC} one needs to 
make a minor technical assumption that the angles $\phi_j$ made between 
the vectors $w_j$ of an admissible sequence and the $y$-axis are bounded 
away from $\pi/2$.  This can be ensured by choosing $\phi_0$ close to the 
$y$-axis, using (\ref{I:Child1}) and the cross product formula to control 
the angles $\angle w_jw_{j+1}$, and then requiring $|w_0|$ large enough.  
\end{Rem}

It would be desirable if the conclusion of Lemma~\ref{L:Good} 
could be strengthened so that the newly constructed vectors satisfy 
(\ref{NotSoGood}) without the ``$-\delta_j$'' in the exponent, 
for then we can use the lemma to construct admissible sequences by 
recursive definition.  However, it can be shown that this stronger 
statement is false.  (This uses a result of Boshernitzan--see the 
appendix to \cite{Ch}.)  Fortunately, the induction can be rescued 
by using a slight variation of the condition (\ref{D:Good}).  

Let $W_j$ be the set of $w\in W$ with the following property: 
for all $t\geq\delta_j$, 
\begin{equation}\label{VeryGood}
  \spec(w)\cap[e^t|w|\log|w|,|w|^{1+t}]\neq\emptyset.  
\end{equation}

\begin{Lem}\label{L:DIO}
  There exists $L_0>0$ such that (\ref{VeryGood}) holds 
  for all $t\geq e_0$ if $|w|\geq L_0$.  
\end{Lem}
\begin{proof}
If (\ref{VeryGood}) does not hold for some $t\geq e_0$, then $w$ has 
convergents $v_k$ and $v_{k+1}$ satisfying $|v_k|<e^t|w|\log|w|$ and 
$|v_{k+1}|>|w|^{1+t}$.  On the one hand we have 
$|w\times v_k|<|w|/|v_{k+1}|<|w|^{-t}$ by (\ref{VCFI}); 
on the other hand we have (\ref{DIO}) implies 
$|w\times v_k|>c_0/|v_k|^{d_0}>c_0e^{-d_0t}|w|^{-d_0}(\log|w|)^{-d_0}$.  
These inequalities contradict each other if 
$$t \geq \frac{d_0\log|w|+d_0\log\log|w|-\log c_0}{\log|w|-d_0}.$$  
Therefore, if $L_0$ is chosen large enough so that the RHS is $<e_0$ 
for $|w|\geq L_0$, then (\ref{VeryGood}) holds for all $t\geq e_0$.  
\end{proof}

\begin{Prop}\label{P:Great}
There exist $L_0>0$ and $\rho_2>0$ such that if $w_j\in W_j$ belongs 
to an admissible sequence $(w_0,\dots,w_j)$ with $|w_0|\geq L_0$, then 
it has $\rho_2|w_j|^{\delta_j}/\log|w_j|$ children contained in the set 
$W_{j+1}$.  
\end{Prop}
\begin{proof}
Let $v_k\in\Spec(w_j)$ be the unique convergent of $w_j$ determined by 
  the condition $|v_k|\leq|w_j|^{1+\delta_j}<|v_{k+1}|$ so that 
$$|v_k|=e^{t_1}|w_j|\log|w_j|\quad\text{and}\quad|v_{k+1}|=|w_j|^{1+t_2}$$
  for some $t_1\geq\delta_j$ and $t_2\leq t_1$.  
If $t_1\geq e_0+\delta_j$ then $w_j$ satisfies the hypothesis of 
Lemma~\ref{L:Good} with $t=e_0+\delta_j$ and each child constructed 
  by the lemma satisfies (\ref{NotSoGood}), which is easily seen to 
  imply (\ref{VeryGood}) for $t\in[\delta_{j+1},e_0]$.  
By Lemma~\ref{L:DIO} it follows that all children constructed lie 
in the set $W_{j+1}$.  Therefore, the conclusion of the proposition 
holds in this case for $\rho_2=\rho_1e^{-2e_0}$.  

Now consider the case $t_1<e_0+\delta_j$.  This time Lemma~\ref{L:Good} 
is applied with $t=t_1$ to obtain the same number of children as before, 
each satisfying (\ref{NotSoGood}) with $t$ replaced by $t_1$.  
Let $W'$ consist of those children which do not belong to $W_{j+1}$.  
Our goal is to show $W'$ occupies only a small fraction (independent of 
$j$) of all the children constructed, provided $L_0$ is large enough.  

Let $\vhi:W'\hookrightarrow Z$ be the function that assigns to any 
$w_{j+1}\in W'$ the unique convergent $v''\in\Spec(w_{j+1})$ with 
maximal Euclidean length $|v''|\leq|w_{j+1}|^{1+\delta_j}$.  The plan 
is to show $\vhi$ is injective then control the cardinality of its image.  

First, we claim every $v''\in\im\vhi$ satisfies an inequality of the form 
\begin{equation}\label{I:v''}
|w_j\times v''\pm ga| \leq 
     \frac{1}{|w_j|^{(1+\delta_j)\max(\delta_{j+1},t_1-\delta_j)}}
\end{equation}
for some positive integer $a<e^{e_0}$ and a choice of the sign on 
the LHS.  Indeed, suppose $v''=\vhi(w_{j+1})$ and let $v'''$ be 
the next convergent after $v''$.  Then     $$|v''|=e^{t_3}|w_{j+1}|
    \log|w_{j+1}|\quad\text{and}\quad|v'''|=|w_{j+1}|^{1+t_4}$$ for 
some real numbers $t_3\geq t_1-\delta_j$, since $w_{j+1}$ satisfies 
(\ref{NotSoGood}), and $t_4>\delta_{j+1}$, by definition of $v'''$.  
Note that $t_4>t_3$ because\footnote{Actually, this requires a short 
calculation because if (\ref{VeryGood}) fails for some $t\geq\delta_{j+1}$, 
then $t_3<t<t_4$ holds \textit{only if} we know 
     $e^t|w_{j+1}|\log|w_{j+1}| \leq |w_{j+1}|^{\delta_{j+1}}$, 
but this holds provided $L_0$ is chosen large enough as required by 
Lemma~\ref{L:eps} for $\eps=e^{-2e_0}$, since $t<e_0+\delta_j\leq2e_0$.}  
$w_{j+1}\in W'$.  By the RHS of (\ref{VCFI}) and the first inequality 
in (\ref{I:Child1}) we have 
$$|w_{j+1}\times v''|\leq\frac{1}{|w_{j+1}|^{t_4}}\leq
         \frac{1}{|w_j|^{(1+\delta_j)\max(\delta_{j+1},t_1-\delta_j)}}.$$ 
Recalling the consecutive pair of convergents $v$ and $v'$ constructed 
in the proof of Lemma~\ref{L:Good} we see that $v''=av'+bv$ for some 
positive integers $a$ and $b$, except in the case when $v''=v'$.  In 
any case $a>0$ and the definition of $v''$ implies these are the only 
possibilities.  Recall also that $v'$ is the convergent responsible 
for (\ref{NotSoGood}) and since $t_1\geq\delta_j$ it follows that 
$|v'|\geq|w_{j+1}|\log|w_{j+1}|$.  On the other hand, we have $t_3<e_0$, 
for otherwise Lemma~\ref{L:DIO} would imply $w_{j+1}\not\in W'$;  
therefore $|v''|<e^{e_0}|w_{j+1}|\log|w_{j+1}|$ and since $|v''|>a|v'|$ 
(because angles between convergents are acute) we have $a<e^{e_0}$.  
Finally, observe that 
$w_{j+1}\times v''=w_j\times v''+gv\times v''=w_j\times v''\pm ga$, 
by (\ref{VPR}).  This proves the claim.  

Next, we show $\vhi$ is injective.  Let $v''_i=\vhi(w_{j+1}^i)$ 
for $i=1,2$ and recall in the proof of Lemma~\ref{L:Good} it was 
shown that $w_{j+1}^i=w_j+gv^i$ for some $v^i\in\Spec(w_{j+1}^i)$.  
Obviously, $v^1\neq v^2$ since we are assuming $w_{j+1}^1\neq w_{j+1}^2$.  
Using (\ref{I:v}), we see that 
$$\angle v^1v^2 \geq \frac{|v^1\times v^2|}{|v^1||v^2|}
                \geq \frac{1}{4|w_j|^{2+2\delta_j}}$$ 
which is approximately a factor of $\log|w_j|$ greater than 
the angle either vector makes with the corresponding child: 
$$\angle v^iw_{j+1}^i = \sin^{-1}\frac{|w_j\times v^i|}{|w_{j+1}^i||v^i|}
                  \gtrapprox \frac{1}{|w_j|^{2+2\delta_j}\log|w_j|}.$$  
Since the angle between $v''_i$ and $w_{j+1}^i$ is even smaller than the 
above, it follows that the vectors $v''_1$ and $v''_2$ are also distinct.  

Finally, we bound the number of vectors in $\im\vhi$.  
Let $v''_i$ for $i=1,2$ be two vectors in $\im\vhi$ satisfying 
(\ref{I:v''}) with the same sign and the same positive integer $a$.  
Put $u=v''_1-v''_2$ and recall the definition of $v_k\in\Spec(w_j)$ 
at the beginning of the proof.  
\textbf{Claim:} If $|u|\leq(1/4)|w_j|^{1+\delta_j}/4$ then $u=\pm dv_k$ 
   for some positive integer $d<4|w_j|^{\delta_j(1-\delta_{j+1})}$.  
The claim implies we either have $$|u|>(1/4)|w_j|^{1+\delta_j} \quad
\text{or}\quad |u|<4e^{2e_0}|w_j|^{1+\delta_j(1-\delta_{j+1})}\log|w_j|.$$ 
Observing that the former is much greater than the latter, which means 
according to (\ref{I:v''}) the vectors in $W'$ are contained in $<2e^{e_0}$ 
narrow strips parallel to $w_j$ and within each strip there are 
$<2(2g+1)e^{2e_0}\log|w_j|$ clusters\footnote{A rough estimate: 
$|u| \leq 2e^{t_3}|w_{j+1}|\log|w_{j+1}| 
       < 2(2g+1)e^{2e_0}|w_j|^{1+\delta_j}\log|w_j|$ provided 
$|w_0|$ is large enough.  Here, we used $w_{j+1}^i\in W'$ and 
Lemma~\ref{L:DIO} to get $|v''_i|<e^{e_0}|w_{j+1}|\log|w_{j+1}|$ 
then apply $|w_{j+1}^i|<(2g+1)|w_j|^{1+\delta_j}$.}  
each having $<4|w_j|^{\delta_j(1-\delta_{j+1})}$ vectors.  
If $L_0$ is chosen large enough as required by Lemma~\ref{L:eps} for 
$\eps=32^{-1}(2g+1)^{-1}e^{-2e_0}\rho_1$, then it follows that less 
than half the children constructed lie in $W'$, i.e. assuming the claim, 
the proposition holds in this case with $\rho_2=(1/2)\rho_1e^{-2e_0}$.  

To prove the claim, we suppose $|u|\leq(1/4)|w_j|^{1+\delta_j}$.  
Let $d=\gcd(u)$ so that $u=dv$ for some $v\in Z$.  
By the triangle inequality and the convenient fact that 
          $(1+\delta_j)\delta_{j+1} > \delta_j$ we have 
\begin{equation}\label{I:wxu}
|w_j\times u| \leq \frac{2}
                   {|w_j|^{(1+\delta_j)\max(\delta_{j+1},t_1-\delta_j)}}
              \leq \frac{2}{|w_j|^{\delta_j}} 
\end{equation}
which together with $|v|\leq|u|\leq(1/4)|w_j|^{1+\delta_j}$ implies 
$$\frac{|w_j\times v|}{|w_j|}\leq\frac{|w_j\times u|}{|w_j|}\leq
\frac{2}{|w_j|^{1+\delta_j}}\leq\frac{1}{2|u|}\leq\frac{1}{2|v|}.$$  
Therefore, $v=\pm v_{k'}$ for some convergent $v_{k'}\in\Spec(w_j)$.  
(See the Remark~\ref{R:VSCC} for an explanation of how the hypothesis 
of Lemma~\ref{VSCC} is satisfied.) 
Since $|v|<|w_j|^{1+\delta_j}$ by hypothesis, $k'\leq k$ by definition 
of $v_k$.  In fact, we must have equality for if $k'<k$ then using the 
LHS of (\ref{VCFI}), the fact that Euclidean lengths of convergents 
form an increasing sequence, and $t_1\leq2e_0$, we have 
$$|w_j\times u| \geq |w_j\times v| \geq
   \frac{|w_j|}{|v_{k'+1}+v_{k'}|} \geq \frac{|w_j|}{2|v_k|} \geq
   \frac{1}{2e^{2e_0}\log|w_j|}$$
which contradicts (\ref{I:wxu}) if $L_0$ is chosen large enough as 
required by Lemma~\ref{L:eps} for $\eps=(1/4)e^{-2e_0}$.  Using 
(\ref{I:wxu}) and the preceding facts about continued fractions 
once again, we get 
$$d = \frac{|w_j\times u|}{|w_j\times v|} < 
 \frac{2|v_{k+1}+v_k|}{|w_j|^{1+(1+\delta_j)\max(\delta_{j+1},t_1-\delta_j)}} 
 < 4|w_j|^{t_2-(1+\delta_j)\max(\delta_{j+1},t_1-\delta_j)}$$
which shows $d<4|w_j|^{\delta_j(1-\delta_{j+1})}$ because 
$$t_2-(1+\delta_j)\max(\delta_{j+1},t_1-\delta_j) \leq 
      t_2-t_1+\delta_j-\delta_j\delta_{j+1} \leq \delta_j(1-\delta_{j+1}).$$  
This proves the claim, and hence the proposition.  
\end{proof}

\begin{proof}[Proof of Theorem~\ref{Thm1}]
Since $\delta_0=e_0$ Lemma~\ref{L:DIO} implies any vector in $W$ 
with large enough Euclidean length belongs to $W_0$.  
Choose any $w_0\in W_0$ with $|w_0|$ greater than the value of 
$L_0$ given by Proposition~\ref{P:Great}.  (In addition, require 
that the initial direction be chosen as in Remark~\ref{R:VSCC}.)  
Applying Proposition~\ref{P:Great} inductively we construct an 
infinite number of admissible sequences $(w_j)$ with the property 
$w_j\in W_j$ for all $j\geq0$.  Moreover, Lemma~\ref{L:Sublinear} 
implies the directions of the vectors in each sequence converge 
to a slowly divergent nonergodic direction with sublinear rate 
$\leq1-1/e_0$.  Any vector $w_j$ occuring at the $j$th stage of 
the construction has at least 
      $$m_j=\rho_2|w_j|^{\delta_j}/\log|w_j|$$ children 
for which the inductive process may be continued indefinitely.  
The angle between the directions of these children are at least 
$$\eps_j=\frac{c}{|w_j|^{2(1+\delta_j)}}$$ where $c>0$ is constant 
depending only on $g$.  Using [Fa, Example~4.6], we see the 
Hausdorff dimension of the set of directions constructed is at least 
\begin{align*}
\liminf_{j\to\infty}\frac{\log m_0\cdots m_{j-1}}{-\log m_j\eps_j} 
   = \liminf_{j\to\infty} 
     \frac{\sum_{i=0}^{j-1}\delta_i\log|w_i|}{(2+\delta_j)\log|w_j|}
   = \lim_{j\to\infty} 
     \frac{1}{2}\left(1-\prod_{i=0}^{j-1}\frac{1}{1+\delta_i}\right)
   = \frac{1}{2}.  
\end{align*}
\end{proof}

\section{Arbitrarily slowly divergent directions}\label{S:Aslow}

To prove Theorem~2 we shall show given any function $R(t)$ with 
$R(t)\to\infty$ as $t\to\infty$ there exists a sequence $(w_j)$ 
satisfying the hypotheses of Propostion~\ref{P:PL} together with 
$M_j\leq R(T_j)$ for $j$ large enough and $\lim m_j=\infty$.  
Note that the length of the shortest simple closed curve on 
$(X_t^\theta,q_t^\theta)$ is at most $2\ell(g_t^\theta V)$ 
since it consists of at most two saddle connections, both 
corresponding to the same vector in $g_t^\theta V$.  Therefore, 
$\lim m_j=\infty$ implies $\lim|w_j|^{-1}w_j$ is divergent.  

The notation $A\asymp B$ means $A/C\leq B\leq AC$ for some implicit 
universal constant $C>0$.  Also, $A\ll B$ means $A\leq B\eps$ for some 
implicit constant $\eps>0$ that may be chosen as small as desired at 
the beginning of the construction.  $B\gg A$ is equivalent to $A\ll B$.  

\begin{Prop}\label{P:Slow}
If $r(t)$ increases to infinity and $r(t)r'(t)\to0$ as $t\to\infty$, 
then there is a sequence $(w_j)_{j\geq0}$ satisfying the hypotheses 
of Proposition~\ref{P:PL} such that for all $j$ large enough 
\begin{enumerate}
\item[(a)] $m_j\leq r(t_j) + C$, where $C:=\log2$, 
\item[(b)] either $|m_{j+1}-r(t_{j+1})|\leq C$ or $m_{j+1}\geq m_j$, 
\item[(c)] if $m_{j+1}<r(t_{j+1})-C$ then $m_{j+1}>m_j+c_0e^{-m_j}$ 
           for some $c_0>0$, and 
\item[(d)] $|w_{j+1}|\asymp|w_j|e^{m_j}$.  
\end{enumerate}
\end{Prop}
\begin{proof}[Proof of Theorem~2 assuming Proposition~\ref{P:Slow}]
Observe the hypotheses of the proposition is satisfied by the 
logarithm of any smooth Lipschitz function increasing to infinity.  
Therefore, given $R(t)$ and any $\eps>0$ one can readily find 
$r(t)$ satisfying the hypotheses and $R(t)\leq(2+\eps)r(t)$.  

Given $r(t)$, (a)-(c) implies $m_j\to\infty$ as $j\to\infty$.  
(d) implies $M_j=m_j+m_{j+1}+O(1)\leq r(t_j)+r(t_{j+1})+O(1)$.  
Using $r(t)r'(t)\to0$ and $t_{j+1}-T_j=m_j+O(1)$ we have 
$r(t_{j+1})\leq r(T_j) + O(1)$ while $r(t_j)\leq r(T_j)$ 
since $r(t)$ is increasing.  Hence, $M_j\leq2r(T_j)+O(1)\leq R(T_j)$.  
\end{proof}

Before proving Proposition~\ref{P:Slow} we need a lemma.  
Note for any pair $(w,v)\in W\times Z$ there exists a unique vector 
$u\in Z$ satisfying $$|w\times u|<\frac{1}{2}
|w\times v|,\quad|u\times v|=1,\quad\text{and}\quad w\cdot u>0.$$ 
Here, we used the fact that vectors in $W$ have irrational slope.  
The next lemma estimates the length of $u$.  
\begin{Lem}\label{L:|u|}
If $|w|>b|v|$ and $|w\times v|\leq\eps\leq1/2\sqrt{2}$ then 
\begin{equation}\label{I:|u|}
  (1-\eps/b)\,|w|\,|w\times v|^{-1} \leq |u| 
       \leq  (1+\eps/b)\,|w|\,|w\times v|^{-1}.  
\end{equation}
\end{Lem}
\begin{proof}
First, consider the case where $w$ lies between $u$ and $v$.  
Then $u+v$ lies between $w$ and $v$ so that comparing the component 
of the vectors $u$, $w$ and $u+v$ orthogonal to $v$ we obtain 
\begin{equation}\label{I:Cfrac}
   \frac{1}{|u+v|} \leq \frac{|w\times v|}{|w|} \leq \frac{1}{|u|}.  
\end{equation}
Let $a,b',c>0$ be given by $|v|=a|u|$, $|v|=b'|w|$ and $c=|w\times v|$.  
Using $|u+v|\leq(1+a)|u|$ we note the LHS implies 
$$|w|\leq |u+v|c \leq (1+a)c|u|$$
But then $|v|\leq (1+a)b'c|u|$ 
so that $|u+v|\leq(1+(1+a)b'c)|u|$.  
Repeating the above argument starting with new estimate on $|u+v|$ 
we get $|u+v|\leq(1+(1+(1+a)b'c))b'c|u|$.  By induction, we find 
$|u+v|\leq|u|/(1-b'c)$.  Since $bb'<1$ this gives the LHS of 
(\ref{I:|u|}) while the RHS holds trivially.  

Now consider the case where $u$ lies between $w$ and $v$.  
Then $w$ lies between $u$ and $u-v$ so that comparing the component 
of $u$ $w$ and $u-v$ orthogonal to $v$ we get 
\begin{equation}
   \frac{1}{|u|} \leq \frac{|w\times v|}{|w|} \leq \frac{1}{|u-v|}.  
\end{equation}
Using $|u-v|\geq(1-a)|u|$ we note the RHS implies 
$$|w|\geq |u-v|c \leq (1-a)c|u|$$
But then $|v|\leq (1-a)b'c|u|$ 
so that $|u-v|\geq(1-(1-a)b'c)|u|$.  
Repeating the above argument starting with new estimate on $|u-v|$ 
we get $|u-v|\geq(1-(1-(1-a)b'c))b'c|u|$.  By induction, we find 
$|u-v|\geq|u|/(1+b'c)$.  Since $bb'<1$ this gives the RHS of 
(\ref{I:|u|}) while the LHS holds trivially.  

There are no other cases because Lemma~\ref{L:!v} implies $|u|>|v|$ 
so that $v$ does not lie between $w$ and $u$.  
\end{proof}

\begin{proof}[Proof of Proposition~\ref{P:Slow}] 
Let $w_0\in W$ be arbitrary.  Since $w_0$ has irrational 
slope, we may choose $w_1=w_0+gv_1$ for some $v_1\in Z$ 
so that $|w_0\times w_1|=g|w_0\times v_1|\ll 1$, and in 
particular, $\leq1/2\sqrt{2}$.  Since $r(t)$ is slowly 
increasing, the choice can be made so that $m_1\gg r(t_1)$.  

Given $(w_j,v_j)\in W\times Z$ with $|w_j\times v_j|<1/2\sqrt{2}$ 
let $u_j$ be the unique vector in $Z$ satisfying 
\begin{equation}\label{D:uj}
  |w_j\times u_j|<\frac{1}{2}|w_j\times v_j|, \quad 
  |u_j\times v_j|=1,\quad\text{and}\quad w_j\cdot u_j>0 
\end{equation}
and define $v^i_{j+1},i=0,1,2$ by 
\begin{equation}\label{D:vi}
  v^1_{j+1} = u_j+\sigma v_j,\qquad
  v^2_{j+1} =2u_j+\sigma v_j \qquad\text{and}\qquad
  v^0_{j+1} = u_j-\sigma v_j
\end{equation}
where $\sigma=+1$ if $w_j$ lies between $u_j$ and $v_j$ and 
$-1$ otherwise.  We note here 
$$\frac{1}{2}|w_j\times v_j|\leq |w_j\times v^1_{j+1}| \leq 
|w_j\times v_j|\leq |w_j\times v^0_{j+1}| \leq 2|w_j\times v_j|$$ 
and $|w_j\times v^2_{j+1}|\leq|w_j\times v^1_{j+1}|$.  

The next pair $(w_{j+1},v_{j+1})$ will be chosen among the three 
possibilities $(w^i_{j+1},v^i_{j+1})$ where $w^i_{j+1}=w_j+gv^i_{j+1}$.  
Note that $u_j$ and $v_j$ are uniquely determined by $w_j$.  
Hence, we may let $\delta=\delta(w_j)\in(0,1/2)$ be defined by 
\begin{equation}\label{D:delta}
  |w_j\times u_j| = \delta|w_j\times v_j|.  
\end{equation}
The index $i\in\{0,1,2\}$ is determined according to the following rule: 
\begin{enumerate}
  \item[(A)] if $m_j>r(t_j)+C$ set $i=0$; 
  \item[(B)] otherwise, choose any $i\in\{0,1\}$ satisfying 
           $|m^i_{j+1}-r(t^i_{j+1})|\leq C$ where 
$$t^i_{j+1}=\frac{1}{2}\log\frac{|w^i_{j+1}|^2}{|w_j\times w^i_{j+1}|}\quad
\text{and}\quad m^i_{j+1}=\frac{1}{2}\log\frac{1}{|w_j\times w^i_{j+1}|}$$
    if possible; if not 
  \item[(C)] let $i=1,2$ be the index realizing the larger of 
             $\delta(w^1_{j+1})$ and $\delta(w^2_{j+1})$.  
\end{enumerate}
The choice made in the case of ambiguity will not matter.  

Note the choice $i=0$ implies $$m_j-C\leq m_{j+1} \leq m_j$$ 
while the choice $i=1$ implies $$m_j\leq m_{j+1} \leq m_j+C.$$  
Similarly, for either choice $i=1,2$, we have $m_{j+1}\geq m_j$.  
If $m_j\leq r(t_j)+C$ for $j=j_0$ then the same holds for 
all $j\geq j_0$.  Since $r(t)$ increases to infinity and 
$m_{j+1}\leq m_j$ whenever (A) is used to choose the next 
vector, $m_j\leq r(t_j)+C$ for some $j$.  This proves (a).  

Since (A) is used to choose the next vector for at most 
finitely many $j$, from some point on the only situation 
when $m_{j+1}<m_j$ is if $i=0$ in (B) is used to choose 
the next vector, but then $|m_{j+1}-r(t_{j+1})|\leq C$.  
This proves (b).  

By choosing $v_1$ so that $r(t_1)\gg1$ we can ensure that 
$m_j\gg1$ for all $j\geq1$, since $r(t)$ is increasing.  
In other words, for any $\eps>0$ we can choose $v_1$ so 
that the sequence of pairs $(w_j,v_j)$ constructed satisfy 
$|w_j\times v_j|\leq\eps$.  We have $|w_j|\in O(\eps|u_j|)$ by 
Lemma~\ref{L:|u|} so that $|v_{j+1}|\asymp|u_j|\asymp|w_j|e^{m_j}$ 
for $\eps$ small enough.  This proves (d) and using 
$|v_{j+1}|\gg|w_j|$ it is readily verified that $(w_j)$ 
satisfies the hypotheses of Proposition~\ref{P:PL}.  

It remains to prove (c).  The hypothesis implies (C) is used 
to choose the next vector with either $i=1,2$.  In this case, 
$$m_{j+1}=m_j+\log1/(1-i\delta(w_j))>m_j+\delta(w_j).$$  
Let $\delta=\delta(w_j)$, $\delta_i=\delta(w^i_{j+1})$ and set 
$\delta'=\max(\delta_1,\delta_2)$.  
We need to show $\delta'\leq\delta$ implies $\delta>c_0e^{-m_j}$.  

Let $\Delta u=u^2_{j+1}-u^1_{j+1}$ where $u^i_{j+1}$ is the unique 
vector in $Z$ satisfying (\ref{D:uj}) with $(w^i_{j+1},v^i_{j+1})$ 
in place of $(w_j,v_j)$.  By definition, we have 
$$|w^i_{j+1}\times u^i_{j+1}| 
  =\delta_i|w_j\times v^i_{j+1}|=\delta_i(1-i\delta)|w_j\times v_j|$$
and $$w^i_{j+1}\times u^i_{j+1} 
          = w_j\times u^i_{j+1} + gv^i_{j+1}\times u^i_{j+1}.$$
Note that $$u^i_{j+1}\times v^i_{j+1} 
   = \sgn(w^i_{j+1}\times v^i_{j+1}) = \sgn(w_j\times v^i_{j+1})$$
does not depend on $i=1,2$ (since $\delta<1/2$).  
Therefore, 
\begin{align}
|w_j\times\Delta u| 
     &= |w^2_{j+1}\times u^2_{j+1} -w^1_{j+1}\times u^1_{j+1}| \\
   &\leq (\delta_2(1-2\delta)+\delta_1(1-\delta))|w_j\times v_j|
     < 2\delta|w_j\times v_j|.\label{I:Du}
\end{align}

Let $\Delta u = du$ where $d=\gcd(\Delta u)$ and $u\in Z$.  
We show $u\not=\pm u_j$ provided $\eps$ was chosen small enough 
at the beginning.  Indeed, using $|v_j|\in O(\eps|u_j|)$ 
\begin{align*}
|u^2_{j+1}| &= (1+O(\eps))~|w^2_{j+1}|~|w_j\times v^2_{j+1}|^{-1}\\ 
  &\geq  (g+O(\eps))~|v^2_{j+1}|~|w_j\times v_j|^{-1}\\ 
  &\geq (2g+O(\eps))~|u_j|~|w_j\times v_j|^{-1} 
\end{align*}
and using $|w_j|\in O(\eps|v^1_{j+1}|)$ and $|v_j|\in O(\eps|u_j|)$ 
\begin{align*}
|u^1_{j+1}| &= (1+O(\eps))~|w^1_{j+1}|~|w_j\times v^1_{j+1}|^{-1}\\ 
  &\leq \frac{g+O(\eps)}{1-2\delta}~    
           |v^1_{j+1}|~|w_j\times v_j|^{-1}\\ 
  &\leq (3g/2+O(\eps))~|u_j|~|w_j\times v_j|^{-1} 
\end{align*}
from which it follows $|\Delta u|>|u_j|$ (provided $\eps$ is small 
enough) so that $d\geq2$.  Hence, $u=\pm u_j$ contradicts (\ref{I:Du}).  

From (\ref{I:Du}) we see that the vector $u'=2u_j+\Delta u$ lies 
between $u_j$ and $w_j$.  Therefore, 
$$\frac{|w_j\times u_j|}{|w_j|} \geq \frac{|u'\times u_j|}{|u'|} 
           \geq \frac{1}{3|\Delta u|}$$
and since $|\Delta u| \in O(|w_j||w_j\times v_j|^{-2})$ 
it follows that $\delta\geq c_0e^{-m_j}$.  
\end{proof}

\end{document}